\input amstex
\documentstyle{amsppt}

\def\QED{\qed}
\expandafter\ifx\csname amsppt.sty\endcsname\relax
\let\QED\relax
\documentstyle{jag}
\issueinfo{00}
{0}
{Xxxx}
{2005}
\fi

\input label.def
\input degt.def
\input cd.def

\def\everyendproof{\QED}

\newtheorem{Remark}\Remark\thm\endAmSdef


{\catcode`\@=11
\gdef\proclaimfont@{\sl}}

\loadbold

\def\todo:{{\bf to do:}}

\let\=\B
\def\ie{\emph{i.e\.}}
\def\eg{\emph{e.g\.}}
\def\cf.{\emph{cf\.}}
\def\all.{\emph{et al\.}}
\def\etc.{\emph{etc\.}}
\def\via{\emph{via}}
\def\dash{\item"\hfill--\hfill"}
\def\Dashes{\widestnumber\item{--}}

\def\PP{\Bbb P}
\def\Aut{\operatorname{Aut}}
\def\discr{\operatorname{discr}}
\def\rank{\operatorname{rk}}
\def\Pic{\operatorname{Pic}}
\def\Br{\operatorname{Br}}
\def\<#1>{\langle#1\rangle}
\def\refl#1{t_{#1}}

\def\dual{^*}

\def\F{\Bbb F}

\let\CL\CalL

\let\CN\CalN

\let\CS\CalS

\let\CU\CalU

\let\CV\CalV

\let\CK\CalK

\let\CH\CalH
\def\gR(#1){\frak r#1}
\let\SS\Sigma
\def\tS{\tilde S}
\def\tSS{\tilde\SS}
\def\tCS{\tilde\CS}
\def\bA{\bold A}
\def\bD{\bold D}
\def\bE{\bold E}
\def\bU{\bold U}
\def\bM{\bold M}
\def\BX{\bar X}
\def\BY{\bar Y}
\def\BC{\bar C}
\def\Bp{\bar p}
\def\Bf{\bar f}
\def\X{_X}
\def\SO{\operatorname{\text{\sl SO}}}
\def\O{\operatorname{\text{\sl O}}}
\def\OAut#1{\O(#1)}
\def\OPlus#1{\O^+(#1)}
\def\Oh#1{\O_h(#1)}
\def\che{\bar\alpha}
\def\cheb#1{\bar\beta^{(#1)}}

\def\barV{\hbox to0pt{$\overline{\phantom{V}}$\hss}V}

\let\Gf\varphi
\let\Go\omega
\let\Gs\sigma
\let\Gr\rho
\let\GO\Omega
\def\GOK{\widetilde\GO}
\def\Cp#1{\Bbb P^{#1}}

\topmatter

\title
On deformations of singular plane sextics
\endtitle

\author
Alex Degtyarev
\endauthor

\address
Department of Mathematics,
Bilkent University,
06800 Ankara, Turkey
\endaddress

\email
degt\@fen.bilkent.edu.tr
\endemail

\abstract
We study complex plane projective sextic curves with simple
singularities up to equisingular deformations. It is shown that
two such curves are deformation equivalent if and only if the
corresponding pairs are diffeomorphic. A way to enumerate all
deformation classes is outlined, and a few examples are
considered, including classical Zariski pairs; in particular,
promising candidates for homeomorphic but not diffeomorphic pairs
are found.
\endabstract



\endtopmatter

\document

\section{Introduction}\label{SS.intro}

\subsection{Motivation and principal results}\label{S.motivation}
Following the real algebraic geometry tradition, an equisingular
deformation of complex plane projective algebraic curves is
called a \emph{rigid isotopy}.
Whenever two curves $C_1,C_2\subset\Cp2$ are rigidly
isotopic, the pairs $(\Cp2,C_i)$, $i=1,2$, are homeomorphic and,
in the case of simple singularities only, also
diffeomorphic.
In his celebrated paper~\cite{Zariski}, O\.~Zariski constructed a
pair of irreducible curves $C_1$, $C_2$ of degree six that have
the same set of singularities (six cusps) but are not rigidly
isotopic; in fact, the complements $\Cp2\sminus C_i$, $i=1,2$, are
not homeomorphic. E\.~Artal~\cite{Artal.def} suggested to call
such curves Zariski pairs. More precisely, a \emph{Zariski pair}
is a pair of reduced plane curves~$C_1$, $C_2$ having the same
combinatorial type of singularities but non-homeomorphic pairs
$(\Cp2,C_i)$; see Section~\ref{S.Zariski.def} for details and
various ramifications. The first degree where
Zariski pairs exist is six, as the rigid isotopy class of a plane
curve of degree up to five is determined by its combinatorial
data, see~\cite{quintics}.

In my thesis (see~\cite{LOMI} and~\cite{poly}), I generalized
Zariski's example and found all pairs of irreducible sextics
$C\subset\Cp2$ that have the same singularities and, as in
Zariski's original case, differ by their Alexander polynomial
(see Section~\ref{S.Zariski.geom} for more
details);
to avoid
confusion with Artal's definition above, we call such curves
\emph{classical Zariski pairs}. I also conjectured that, up to
equisingular deformation, an irreducible sextic is determined by
its set of singularities and its Alexander polynomial. (The
conjecture was based on the calculation for a few special cases
and the fact that the assertion
does hold if the curves have at least one non-simple singular
point, see~\cite{quartics}.)
The conjecture was soon disproved by
H\.~Tokunaga~\cite{Tokunaga}, who constructed a pair of irreducible
sextics~$C_1$, $C_2$ with the same sets of singularities and
Alexander polynomials. Still,
Tokunaga's curves differ by the fundamental group
$\pi_1(\Cp2\sminus C_i)$. In a recent series of
papers~\cite{Artal.first}--\cite{Artal.last} Artal \all.
constructed a number of new examples of not rigidly isotopic
pairs~$(C_1,C_2)$
of sextics; for many
pairs the fundamental groups
$\pi_1(\Cp2\sminus C_i)$ are calculated and shown to coincide.
Thus, the question arises whether the curves
constitute Zariski pairs, \ie, whether $(\Cp2,C_1)$
and $(\Cp2,C_2)$ are homeomorphic. We show that they are not
diffeomorphic. More precisely, the following theorem holds.

\theorem\label{th.main}
Two sextic curves $C_1,C_2\subset\Cp2$ with simple singularities
only are rigidly isotopic if and only if there exists a
diffeomorphism $f\:(\Cp2,C_1)\to(\Cp2,C_2)$ that is \emph{regular}
in the sense that each singular point of~$C_1$ has a
neighborhood~$U$ such that the restriction~$f|_U$ is complex
analytic.
\endtheorem

This theorem is proved in Section~\ref{main.proof}.

\Remark\label{rem.diffeo}
The requirement that $f$ should be a diffeomorphism is not a mere
technical assumption; it is used essentially in the proof as a
means of comparing the orientations of the homological types
of~$C_1$ and~$C_2$ (see Section~\ref{S.homology}).
Since pairs of sextics that differ solely by the
orientation of their
homological types do exist (\eg, Proposition~\ref{A18+A1}), one may anticipate that they would
provide examples of homeomorphic but not diffeomorphic pairs.
\endRemark

As Theorem~\ref{th.main} settles the relative
$\text{Dif}=\text{Def}$ problem for plane sextics, it simplifies
the process of finding Zariski pairs. For example, according to
J.-G\.~Yang~\cite{Yang} there is a five page long list of sextics
with maximal total Milnor number $\mu=19$. The rigid isotopy classes of
such curves are described by definite lattices, which tend to have
very few isometries; hence, there should be a great deal of not
rigidly isotopic pairs sharing the same sets of singularities.

The proof of Theorem~\ref{th.main} is based on an explicit
description of the moduli space of sextics, see
Theorems~\ref{2.components} and~\ref{enumeration}, which, in turn,
is a rather standard application of the global Torelli theorem for
$K3$-surfaces and the surjectivity of the period map. As another
application, Theorem~\ref{enumeration} reduces the rigid isotopy
classification of plane sextics to an arithmetic question about
lattices. We outline the principal steps of enumerating abstract
homological types, see Section~\ref{S.arithmetics}, and apply the
scheme to two polar cases, those of curves with few singularities
and curves with many singularities. In the former case, we prove
Corollary~\ref{cor.l+mu<=19} and Theorem~\ref{th.l+mu<=19}, which
give simple sufficient conditions for a set of
singularities/configuration to be realized by a single rigid
isotopy class. As a further application, we enumerate all curves
constituting classical Zariski pairs without nodes, see
Theorem~\ref{th.Zariski}. For Zariski's original example the
theorem states that {\proclaimfont plane sextics with six cusps
form exactly two deformation families.} To my knowledge, this fact
is new: contrary to the common belief,
Zariski himself has only asserted the existence of \emph{at least}
two families.

In the latter case (maximal total Milnor number $\mu=19$), the
problem reduces to enumerating certain positive definite lattices
of rank~$2$
and their isometries.
The algorithm can easily be implemented
(in fact, I do have it implemented in {\tt Maple}), and, when
combined with Yang's algorithm~\cite{Yang} for enumerating
the configurations, it should produce
a complete list of rigid
isotopy classes.
However, instead of
compiling a long computer aided table, I illustrate the approach
by studying a few examples (see
Propositions~\ref{examples.first}--\ref{examples.last}) that were
first considered in~\cite{Artal.first}--\cite{Artal.last}.

Undoubtedly, the most remarkable example is that given by
Proposition \ref{A18+A1}, where two curves
differ by the orientation of
their homological types.
It is worth mentioning
that found in the literature are a great number of various
deformation classification problems related to the global Torelli
theorem for $K3$-surfaces (in the real case, see recent
papers~\cite{Nik05} and~\cite{finite} and the survey~\cite{survey}
for further references; in the complex case, see, \eg,
V\.~Nikulin~\cite{NikK3}, A\.~Degtyarev \all.~\cite{finite},
Sh\.~Mukai~\cite{Mukai}, Sh\.~Kond\=o~\cite{Kondo1}
and~\cite{Kondo2}, and G\.~ Xiao~\cite{Xiao}). To my knowledge,
the study of singular plane sextics is the only case so far where
the orientation of maximal positive definite subspaces is involved
in an essential way!

\subsection{Contents of the paper}\label{S.contents}
In \S\ref{SS.lattices}, we outline the principal notions and
results of Nikulin's theory of discriminant forms of even
integral lattices. It is largely based on Nikulin's original
paper~\cite{Nikulin}. A preliminary calculation involving certain
definite lattices is also made here. In \S\ref{SS.moduli}, the
relation between plane sextics and $K3$-surfaces is explained, the
moduli space is described, and Theorem~\ref{th.main} is proved.
In \S\ref{SS.sextics}, we discuss a few results relating the
geometry of a sextic and the arithmetic properties of its
homological type.
Finally, \S\ref{SS.examples} deals with the classification of
oriented abstract homological types, which enumerate the rigid
isotopy classes of sextics. We outline the general scheme and
apply it to a few particular examples.

\subsection{Acknowledgements}\label{S.acknowledgements}
I am thankful to S\.~Orevkov, who drew my attention to the
problem, and to E\.~Artal, who introduced me to the modern state
of the subject and encouraged me to develop and publish the results.
My special gratitude is to A\.~Klyachko for his
patient explanation of the $p$-adic machinery behind
Nikulin's
results on the discriminant forms of even
integral lattices. This paper originated from a sample calculation
used in my lecture at the \emph{P\'eriode sp\'eciale de DEA
``Topologie des vari\'et\'es alg\'ebriques r\'eelles''} at
\emph{Universit\'e Louis Pasteur}, Strasbourg.
I am grateful to the organizers of this event for their
hospitality and to the audience for their patience.

\section{Integral lattices}\label{SS.lattices}

\subsection{Finite quadratic forms}\label{S.finite}
A \emph{finite quadratic form} is a finite abelian
group~$\CL$ equipped with a nonsingular quadratic form, \ie, a
map
$q\:\CL\to\Q/2\Z$
satisfying $q(x+y)=q(x)+q(y)+2b(x,y)$
for all $x,y\in\CL$ and some nonsingular symmetric bilinear form
$b\:\CL\otimes\CL\to\Q/\Z$.
If $q$ is understood, we write~$x^2$ and
$x\cdot y$ for $q(x)$ and $b(x,y)$, respectively.

The bilinear form~$b$ is determined by~$q$; it is called the
bilinear form \emph{associated} with~$q$, and $q$ is called a
\emph{quadratic extension} of~$b$.

The group of automorphisms of~$\CL$ preserving~$q$ is denoted
by~$\Aut\CL$.

The \emph{Brown invariant} of a finite quadratic form~$\CL$ is the
residue $\Br\CL\in\Z/8\Z$ defined \via\ the Gauss sum
$$
\exp\bigl(\tfrac14i\pi\Br\CL\bigr)=
 \mathopen|\CL\mathclose|^{-\frac12}\sum_{x\in\CL}\exp\bigl(i\pi x^2\bigr).
$$
The Brown invariant is additive:
$\Br(\CL_1\oplus\CL_2)=\Br\CL_1+\Br\CL_2$.

Clearly, each finite quadratic form~$\CL$ splits canonically into
orthogonal sum of its primary components:
$\CL=\bigoplus\CL\otimes\Z_p$, summation over all primes~$p$.
For a prime~$p$, let $\CL_p=\CL\otimes\Z_p$
be the $p$-primary part of~$\CL$.
Denote by $\ell(\CL)$ the minimal number of generators of~$\CL$,
and let $\ell_p(\CL)=\ell(\CL_p)$. Obviously,
$\ell(\CL)=\max_p\ell_p(\CL)$.

For a fraction $\frac mn\in\Q/2\Z$ with $(m,n)=1$ and
$mn=0\bmod2$, let~$\<\frac mn>$ be the nondegenerate
quadratic form on~$\Z/n\Z$ sending the generator to~$\frac mn$.
For an integer $k\ge1$, let $\CU_{2^k}$ and $\CV_{2^k}$ be the
quadratic forms on the group $(\Z/2^k\Z)^2$ defined by the
matrices
$$
\CU_{2^k}=
 \bmatrix0&\alpha_k\\\alpha_k&0\endbmatrix,\qquad
\CV_{2^k}=
 \bmatrix\alpha_{k-1}&\alpha_k\\\alpha_k&\alpha_{k-1}\endbmatrix,\quad
 \text{where $\alpha_k=\dfrac1{2^k}$}.
$$
(When speaking about the matrix of a finite quadratic form, we
assume that the diagonal elements are defined modulo~$2\Z$ whereas
all other elements are defined modulo~$\Z$.) According to
Nikulin~\cite{Nikulin}, each finite quadratic form is an
orthogonal sum of cyclic summands $\<\frac mn>$ and summands of
the form $\CU_{2^k}$, $\CV_{2^k}$. The Brown invariants of these
elementary blocks are as follows: if $p$ is an odd prime,
then
$$
\Br\Bigl<\frac{2a}{p^{2s-1}}\Bigr>=
 2\Bigl(\frac{a}p\Bigr)-\Bigl(\frac{-1}p\Bigr)-1,\quad
\Br\Bigl<\frac{2a}{p^{2s}}\Bigr>=0\quad
 \text{(for $s\ge1$ and $(a,p)=1$)}.
$$
If $p=2$, then
$$
\gather
\Br\Bigl<\frac{a}{2^k}\Bigr>=a+\frac12k(a^2-1)\bmod8\quad
 \text{(for $k\ge1$ and odd $a\in\Z$)},\\\noalign{\smallskip}
\Br\CU_{2^k}=0,\quad\Br\CV_{2^k}=4k\bmod8\quad
 \text{(for all $k\ge1$)}.
\endgather
$$
Quite a number of relations, \ie, isomorphisms between various
combinations of the aforementioned forms, is also listed
in~\cite{Nikulin}.
These observations make the classification of finite quadratic
forms rather straightforward, although tedious. Two simple known
results used in the sequel are listed below. More details on
quadratic forms on $2$-primary groups can be found in~\cite{DIK}
and~\cite{Nik86}.

\proposition\label{forms.Zp}
Let $p\ne2$ be an odd prime. Then
a quadratic form on a group~$\CL$
of exponent~$p$
is determined by its rank $\ell(\CL)=\ell_p(\CL)$
and Brown invariant $\Br\CL$.
\endproposition

A finite quadratic form
is called \emph{even} if $x^2$ is an integer for each element
$x\in\CL$ of order~$2$; otherwise, it is called \emph{odd}.
Clearly, a form is odd if and only if it contains $\<\pm\frac12>$
as an orthogonal summand.

\proposition[Proposition \rm(see~\cite{Wall} or \cite{Guillou-Marin})]\label{forms.Z2}
A quadratic form on a group~$\CL$
of exponent~$2$
is determined by its rank $\ell(\CL)=\ell_2(\CL)$, parity
\rom(even or odd\rom), and Brown invariant $\Br\CL$.
\endproposition

\subsection{Even integral lattices and discriminant forms}\label{S.lattices}
An \emph{\rom(integral\rom) lattice} is a finitely generated free
abelian group~$L$
equipped with a symmetric bilinear form $\Gf\:L\otimes L\to\Z$.
When the form is understood, we
will freely use the multiplicative notation $u\cdot v=\Gf(u,v)$ and
$u^2=\Gf(u,u)$. A lattice~$L$ is called \emph{even} if
$u^2=0\bmod2$ for each $u\in L$; otherwise, it is called
\emph{odd}.

Since the transition matrix from one integral basis to another one
has determinant~$\pm1$, the \emph{determinant} $\det
L=\det\Gf\in\Z$ is well defined. The lattice~$L$ is called
\emph{non-degenerate} if $\det L\ne0$; it is called
\emph{unimodular} if $\det L=\pm1$. The \emph{signature} of a
non-degenerate lattice~$L$ is the pair $(\Gs_+L,\Gs_-L)$ of its
inertia indices. Recall that $\Gs_+L$ is the dimension of any
maximal positive definite subspace of the vector space
$L\otimes\R$. Recall, further, that all maximal positive definite
subspaces of $L\otimes\R$ can be oriented in a coherent way. For
example, the orientations of two such subspaces $\Go_1$, $\Go_2$
can be compared using the orthogonal projection $\Go_2\to\Go_1$,
which is necessarily
injective and hence bijective.

Given a lattice~$L$, we denote by $\OAut L$ the group of
isometries of~$L$, and by $\OPlus L\subset\OAut L$
its subgroup consisting of the
isometries preserving the
orientation of maximal positive definite subspaces. Clearly,
either $\OPlus L=\OAut L$ or $\OPlus L\subset\OAut L$ is a subgroup of
index~$2$. In the latter case, each element of
$\OAut L\sminus\OPlus L$ is called a \emph{$+$-disorienting} isometry.
(The awkward terminology is chosen to avoid
confusion with
isometries reversing the orientation of $L$
itself.)

If $L$ is a non-degenerate lattice,  the dual group
$L\dual=\Hom(L,\Z)$ can be identified with the subgroup
$\bigl\{x\in L\otimes\Q\bigm|x\cdot y\in\Z\text{ for all $y\in
L$}\bigr\}$. The quotient $L\dual/L$ is called the {\it
discriminant group\/} of~$L$ and is denoted by~$\CL$ or $\discr L$.
One has
$\mathopen|\CL\mathclose|=\mathopen|\det L\mathclose|$ and
$\ell(\CL)\le\rank L$.
The discriminant group inherits from $L\otimes\Q$
a non-degenerate symmetric bilinear
form $b\:\CL\otimes\CL\to\Q/\Z$ and, if $L$ is even, its
quadratic extension $q\:\CL\to\Q/2\Z$. Thus, the discriminant of
an even lattice is a finite quadratic form.

Two integral lattices~$L_1$, $L_2$ are said to have the same
\emph{genus} if all their localizations $L_i\otimes\R$ and
$L_i\otimes\Q_p$ are isomorphic (over~$\R$ and~$\Q_p$,
respectively). Each genus is known to contain finitely many
isomorphism classes. The relation between the genus of a lattice
and its discriminant form is given by the following two
statements (see also Section~\ref{S.existence} below).

\theorem[Theorem \rm(see~\cite{Nikulin})]
The genus of an even integral lattice~$L$ is determined by its
signature $(\Gs_+L,\Gs_-L)$ and discriminant form $\discr L$.
\endtheorem

In what follows, the genus of even integral lattices determined by
a signature $(\Gs_+,\Gs_-)$ and a discriminant form~$\CL$ is
referred to as the \emph{genus $(\Gs_+,\Gs_-;\CL)$}.

\theorem[Theorem \rm(van der Blij formula, see~\cite{vdBlij})]
For any nondegenerate even integral lattice~$L$ one
has $\Br\CL=\Gs_+L-\Gs_-L\bmod8$.
\endtheorem

Since the construction of the discriminant
form~$\CL$ is natural, there is a canonical homomorphism
$\OAut L\to\Aut\CL$. Its image is denoted by $\Aut_L\CL$. Of
special importance are so called reflections of~$L$: given a vector
$a\in L$, the \emph{reflection} against the hyperplane orthogonal
to~$a$ (for short, reflection defined by~$a$) is the automorphism
$$
\refl{a}\:L\to L,\quad
x\mapsto x-2\frac{a\cdot x}{a^2}a.
$$
It is easy to see that $\refl{a}$ is an involution, \ie,
$\refl{a}^2=\id$. The reflection~$\refl{a}$ is well defined
whenever $a\in(a^2/2)L\dual$. In particular, $\refl{a}$ is well
defined if $a^2=\pm1$ or~$\pm2$; in this case the induced
automorphism of the discriminant group~$\CL$ is the identity and
$\refl{a}$ extends to any lattice containing~$L$.

\subsection{Special lattices and notation}\label{S.notation}
Given a lattice~$L$ and an integer~$n$, we denote by~$L(n)$ the
lattice obtained by multiplying all values by~$n$ (\ie, the
quadratic form $x\mapsto nx^2$ defined on the same group~$L$). For
finite quadratic forms the multiplication operation is meaningful
only for $n=-1$, and we abbreviate $-\CL=\CL(-1)$.

The notation~$nL$, $n\ge1$, stands for the direct sum of
$n$~copies of~$L$.

The \emph{hyperbolic plane} is the lattice~$\bU$ spanned by two
vectors~$u$, $v$ so that $u^2=v^2=0$, $u\cdot v=1$. Any pair
$(u,v)$ as above is called a \emph{standard basis} for~$\bU$. In
fact, it is unique up to transposing~$u$ and~$v$ and
multiplying one or both of them by~$(-1)$. The hyperbolic plane
is an even unimodular lattice of signature $(1,1)$.

A \emph{root} in a lattice~$L$ is an element $v\in L$ of
square~$-2$.
Given~$L$, we denote by $\gR(L)\subset L$ the sublattice
generated by all roots of~$L$.
A \emph{root system} is a negative definite lattice
generated by its roots. Every root system admits a unique decomposition
into an orthogonal sum of irreducible root systems, the latter
being either $\bA_p$, $p\ge1$, or~$\bD_q$, $q\ge4$,
or~$\bE_6$, $\bE_7$, $\bE_8$. The discriminant forms are as
follows:
$$
\gather
\discr\bA_p=\<-\tfrac{p}{p+1}>,\quad
\discr\bD_{2k+1}=\<-\tfrac{2k+1}4>,\\
\discr\bD_{8k\pm2}=2\<\mp\tfrac12>,\quad
\discr\bD_{8k}=\CU_2,\quad
\discr\bD_{8k+4}=\CV_2,\quad\\
\discr\bE_6=\<\tfrac23>,\quad
\discr\bE_7=\<\tfrac12>,\quad
\discr\bE_8=0.
\endgather
$$
The orthogonal group of a root system~$L$ is a semi-direct
product of the group generated by reflections (defined by the
roots of~$L$), which acts simply transitively on the set of Weyl chambers
of~$L$, and the group of symmetries of any fixed Weyl chamber (or
Dynkin graph) of~$L$. As a consequence, the following
statement holds:

\proposition\label{root.auto}
For a root system~$L$,
the subgroup $\Aut_L\CL$  coincides with the image in
$\Aut\CL$ of the group of symmetries of any fixed Weyl chamber.
\endproposition

If $L$ is an irreducible root systems other than~$\bA_p$
or~$\bD_q$ with $q=8k+4\ge12$, one has $\Aut_L\CL=\Aut\CL$. If
$L=\bA_p$, the image $\Aut_L\CL$ is the subgroup $\{\pm\id\}$. In
the case $L=\bD_{8k+4}$, $k\ge1$, the full orthogonal group
$\Aut\CL$ is the group~$S_3$ of permutations of the three elements
of square~$1\bmod2\Z$, whereas the image $\Aut_L\CL$ is generated
by one of the three transpositions.

Further details on irreducible root systems are found in
N\.~Bourbaki~\cite{Bourbaki}.

\subsection{Definite lattices of rank~$2$}\label{S.rank2}
Each positive definite even lattice~$N$ of rank~$2$ has a unique
representation by a matrix of the form
$$
\bmatrix2a&b\\b&2c\endbmatrix,\quad
0<a\le c,\quad0\le b\le a.\eqtag\label{matrix.2}
$$
Denote the lattice represented by~\eqref{matrix.2} by
$\bM(a,b,c)$.
Let $(u,v)$ be a basis in which the quadratic form is given
by~\eqref{matrix.2}. Then, depending on~$a$, $b$, and~$c$,
the orthogonal group $\OAut N$ is one of the groups described below:
\Dashes
\roster
\dash
\hbox{$0<b<a<c$: }the group $\OAut N\cong\Z/2\Z$ is generated by~$-\id$;

\dash
\hbox{$0<b<a=c$: }the group $\OAut N\cong\Z/2\Z\times\Z/2\Z$ is
generated by~$-\id$ and the transposition $(u,v)\mapsto(v,u)$;

\dash
\hbox{$b=0$, $a<c$: }the group $\OAut N\cong\Z/2\Z\times\Z/2\Z$ is
generated by~$\refl{u}$ and~$\refl{v}$;

\dash
\hbox{$b=0$, $a=c$: }then $N=2\bA_1(-a)$ and $\OAut N\cong\Bbb D_4$ is the
group of symmetries of a square; it is generated by~$\refl{u}$ and
the transposition $(u,v)\mapsto(v,u)$;

\dash
\hbox{$b=a<c$: }the group $\OAut N\cong\Z/2\Z\times\Z/2\Z$ is
generated by~$-\id$ and~$\refl{u}$;

\dash
\hbox{$b=a=c$: }then $N=\bA_2(-a)$ and $\OAut N\cong\Bbb D_6$ is the
group of symmetries of a regular hexagon; it is generated
by~$\refl{u}$ and the transposition
$(u,v)\mapsto(v,u)$.

\endroster

All results above are classical and well known. The inequalities
$a\le c$ and $\mathopen|b\mathclose|\le a$ can be achieved by a
sequence of transpositions $(u,v)\mapsto(v,u)$ and transformations
$(u,v)\mapsto(u,v\pm u)$. Then, assuming that the matrix has the
form~\eqref{matrix.2}, for a vector $xu+yv\in N$ one has
$$
(xu+yv)^2=2ax^2+2bxy+2cy^2\ge
 2a(x^2+y^2)-2a\mathopen|xy\mathclose|\ge
 a(x^2+y^2).
$$
Since $x$ and~$y$ are integers, it immediately follows that $u$ is
a shortest vector and, unless $a=c$, the only shortest vectors are
$\pm u$. If $a=c$, there are two more shortest vectors~$\pm v$,
and if also $b=a$, there are yet two more, $\pm(u-v)$. From here,
one can easily deduce the uniqueness of
representation~\eqref{matrix.2}. The description of the
orthogonal group is also straightforward: one observes that~$u$
should be taken to a shortest vector and then, assuming~$u$ fixed,
the only nontrivial isometry of the Euclidean plane $N\otimes\R$
is the reflection against the line spanned by~$u$; it remains to
enumerate the few cases when this reflection is defined over~$\Z$.

\subsection{Nikulin's existence and uniqueness results}\label{S.existence}
Let $p$ be a prime. The notion of lattice and its discriminant
form extends to the case of finitely generated free
$\Z_p$-modules. (In the case $p=2$, to define the quadratic form
on the discriminant group one still needs to require that the
lattice should be even.) The discriminant of a $p$-adic
lattice~$L_p$ is a finite $\Z_p$-module~$\CL_p$ (in other words,
$p^k\CL_p=0$ for some~$k$ large enough), and one has
$\mathopen|\CL_p\mathclose|=\mathopen|\det
L_p\mathclose|\bmod\Z_p^*$. For an integral lattice~$L$ one has
$\discr(L\otimes\Z_p)=(\discr L)\otimes\Z_p=\CL_p$.

According to Nikulin~\cite{Nikulin}, given a prime~$p$ and a
$\Q/2\Z$-valued quadratic form on a finite $\Z_p$-module~$\CL$,
there is a $p$-adic lattice $L$ such that $\rank L=\ell_p(\CL)$
and $\discr L=\CL$. Unless $p=2$ and $\CL$ is odd, such a
lattice~$L$ is determined by~$\CL$ uniquely up to isomorphism; in
particular, the ratio $\det L/\mathopen|\CL\mathclose|$ is a well
defined element of the group $\Z_p^*/(\Z_p^*)^2$. We will denote
it by $\det_p\CL$.
In the exceptional case $p=2$, $\CL$~odd there are
two lattices~$L$ as above, the ratio of their determinants being
$5\in\Z_2^*/(\Z_2^*)^2$.

\theorem[Theorem \rm(see Theorem~1.10.1
in~\cite{Nikulin})]\label{N.existence}
Let $\CL$ be a finite quadratic form and let $\Gs_\pm$ be a
pair of integers.
Then, the following four conditions are necessary and sufficient
for the existence of an even integral lattice~$L$ whose
signature is $(\Gs_+,\Gs_-)$ and whose discriminant form
is~$\CL$\rom:
\roster
\item\local{>0}
$\Gs_\pm\ge0$ and $\Gs_++\Gs_-\ge\ell(\CL)$\rom;
\item\local{Br}
$\Gs_+-\Gs_-=\Br\CL\bmod8$\rom;
\item\local{p}
for each $p\ne2$, either $\Gs_++\Gs_->\ell_p(\CL)$ or
$\det_p\CL_p=(-1)^{\Gs_-}\bmod(\Z_p^*)^2$\rom;
\item\local2
either $\Gs_++\Gs_->\ell_2(\CL)$, or
$\CL_2$ is odd, or
$\det_2\CL_2=\pm1\bmod(\Z_2^*)^2$.
\endroster
\endtheorem

\theorem[Theorem \rm(see Theorem~1.13.2
in~\cite{Nikulin})]\label{N.uniqueness}
Let~$L$ be an indefinite even integral lattice, $\rank L\ge3$. The
following two conditions are sufficient for~$L$ to be unique in
its genus\rom:
\roster
\item\local{p}
for each $p\ne2$, either $\rank L\ge\ell_p(\CL)+2$ or
$\CL_p$ contains a subform isomorphic to $\<a/p^k>\oplus\<b/p^k>$,
$k\ge1$, as an orthogonal summand\rom;
\item\local2
either $\rank L\ge\ell_2(\CL)+2$ or
$\CL_2$ contains a subform isomorphic to $\CU_{2^k}$, $\CV_{2^k}$, or
$\<a/2^k>\oplus\<b/2^{k+1}>$, $k\ge1$, as an orthogonal summand.
\endroster
\endtheorem

\theorem[Theorem \rm(see Theorem~1.14.2
in~\cite{Nikulin})]\label{N.strong.uniqueness}
Let~$L$ be an indefinite even integral lattice, $\rank L\ge3$. The
following two conditions are sufficient for~$L$ to be unique in
its genus and for the canonical
homomorphism $\OAut L\to\Aut\CL$ to be onto\rom:
\roster
\item\local{p}
for each $p\ne2$, $\rank L\ge\ell_p(\CL)+2$\rom;
\item\local2
either $\rank L\ge\ell_2(\CL)+2$ or
$\CL_2$ contains a subform isomorphic to $\CU_2$ or $\CV_2$
as an orthogonal summand.
\endroster
\endtheorem

\subsection{Extensions}\label{S.extension}
 From now on we confine ourselves to even lattices.
An \emph{extension} of an even lattice~$S$ is an even lattice~$L$
containing~$S$. Two extensions $L_1\supset S$ and $L_2\supset S$
are called \emph{isomorphic} (\emph{strictly isomorphic}) if there
is an isomorphism $L_1\to L_2$ preserving~$S$ (respectively,
identical on~$S$). More generally, one can fix a subgroup
$A\subset\OAut S$ and speak about \emph{$A$-isomorphisms} and
\emph{$A$-automorphisms} of extension,
\ie, isometries
whose restriction to~$S$ belongs to~$A$.

Any extension $L\supset S$ of finite index admits a unique
embedding $L\subset S\otimes\Q$. If $S$ is nondegenerate, then $L$
belongs to~$S\dual$ and thus defines a subgroup
$\CK=L/S\subset\CS$, called the \emph{kernel} of the extension.
Since $L$ itself is an integral lattice, the kernel~$\CK$ is
isotropic, \ie, the restriction to~$\CK$ of the discriminant
quadratic form is identically zero. Conversely, given an isotropic
subgroup $\CK\subset\CS$, the subgroup
$$
L=\bigl\{x\in S\dual\mid(x\bmod S)\in\CK\bigr\}\subset S\dual
$$
is an extension of~$S$.
Thus, the following statement holds:

\proposition[Proposition \rm(see~\cite{Nikulin})]\label{extension.finite}
Let~$S$ be a nondegenerate even lattice, and fix a subgroup
$A\subset\OAut S$.
The map $L\mapsto\CK=L/S\subset\CS$ establishes a
one-to-one correspondence between the set of $A$-isomorphism
classes of finite index extensions $L\supset S$ and the set
of $A$-orbits
of isotropic subgroups $\CK\subset\CS$. Under this correspondence,
one has $\discr L=\CK^\perp/\CK$.

An isometry $f\:S\to S$ extends to a finite index extension
$L\supset S$ defined by an isotropic subgroup $\CK\subset\CS$ if
and only if the automorphism $\CS\to\CS$ induced by~$f$
preserves~$\CK$ \rom(as a set\rom).
\endproposition

\Remark
Since a finite index extension $L\supset S$
has the same
signature as~$S$, Proposition~\ref{extension.finite} implies, in
particular, that $\Br(\CK^\perp/\CK)=\Br\CS$ for any isotropic
subgroup $\CK\subset\CS$. This observation facilitates the
calculation of the Brown invariant; for example, it can be used to
reduce the list of values of~$\Br$ given in Section~\ref{S.finite}
to a few special cases.
\endRemark

\corollary\label{extension.roots}
Any imprimitive extension of 
a root system
$S=3\bA_2$, $\bA_5\oplus\bA_2$, $\bA_8$, $\bE_6\oplus\bA_2$,
$2\bA_4$, $\bA_5\oplus\bA_1$, $\bA_7$, $\bD_8$, $\bE_7\oplus\bA_1$, $4\bA_1$,
$\bA_3\oplus2\bA_1$, or
$\bD_q\oplus2\bA_1$ with $q<12$ or $q\ne0\bmod4$ contains a finite
index extension $R\supset S$, where $R$ is a root system strictly
larger than~$S$.
\endcorollary

\proof
The extensions are easily enumerated using
Proposition~\ref{extension.finite}. (In fact, in all cases except
$S=\bD_8\oplus2\bA_1$, a nontrivial finite order extension is
unique up to isometry.) The statement follows then from a
direct calculation, using the fact that each lattice~$\bE_6$,
$\bE_7$, $\bE_8$ is unique in its genus and the known embedding
$2\bA_1\subset\bD_q$ with $(2\bA_1)^\perp_{\bD_q}\cong\bD_{q-2}$.
\endproof

Another extreme case is when $S\subset L$ is a primitive
nondegenerate sublattice and $L$ is a unimodular lattice. Then $L$
is a finite index extension of the orthogonal sum $S\oplus
S^\perp$, both~$S$ and~$S^\perp$ being primitive in~$L$. Since
$\discr L=0$, the kernel $\CK\subset\CS\oplus\discr S^\perp$
is
the graph of an anti-isometry $\CS\to\discr S^\perp$. Conversely,
given a lattice~$N$ and an anti-isometry $\kappa\:\CS\to\CN$, the
graph of~$\kappa$ is an isotropic subgroup
$\CK\subset\CS\oplus\CN$ and the resulting extension $L\supset
S\oplus N\supset S$ is a unimodular primitive extension of~$S$
with $S^\perp\cong N$.

Let~$N$ and $\kappa\:\CS\to\CN$ be as above, and let $s\:S\to S$
and $t\:N\to N$ be a pair of isometries. Then the direct sum
$s\oplus t\:S\oplus N\to S\oplus N$ preserves the graph
of~$\kappa$ (and, thus, extends to~$L$) if and only if
$\kappa\circ s=t\circ\kappa$. (We use the same notation~$s$ and~$t$
for the induced homomorphisms on~$\CS$ and~$\CN$, respectively.)
Summarizing, one obtains the following statement:

\proposition[Proposition \rm(see~\cite{Nikulin})]\label{extension.unimodular}
Let~$S$ be a nondegenerate even lattice, and let $s_+$, $s_-$ be
nonnegative integers. Fix a subgroup $A\subset\OAut S$. Then the
$A$-isomorphism class of a primitive extension $L\supset S$ of~$S$
to a unimodular lattice~$L$ of signature~$(s_+,s_-)$ is determined
by
\roster
\item\local1
a choice of a lattice~$N$
in the genus
$(s_+-\Gs_+S,s_--\Gs_-S;-\CS)$, and

\item\local2
a choice of a bi-coset of the canonical left-right action of
$A\times\Aut_N\CN$ on the set of anti-isometries $\CS\to\CN$.

\endroster
If a lattice~$N$ and an anti-isometry $\kappa\:\CS\to\CN$ as above
are chosen \rom(and thus an extension~$L$ is fixed\rom), an isometry
$t\:N\to N$ extends to an $A$-automorphism of~$L$ if and only if
the composition $\kappa^{-1}\circ t\circ\kappa\in\Aut\CS$ is in
the image of~$A$.
\endproposition

\Remark\label{rem.M-V}
Proposition~\ref{extension.unimodular} can be regarded as the algebraic
counterpart of the Meyer-\hskip0ptVietoris exact sequence of the gluing of
two $4$-manifolds \via\ a diffeomorphism of their boundaries. The
lattices in question are the intersection index forms on the
$2$-homology of the manifolds, and the discriminant forms are the
linking coefficient forms on the $1$-homology of the boundary. The
anti-isometry~$\kappa$ as above is the homomorphism induced by the
identification of the boundaries (which is orientation reversing).
For more details, see, \eg,
O\.~Ivanov and N\.~Netsvetaev~\cite{Nikita1} and~\cite{Nikita2}.
\endRemark

\section{The moduli space}\label{SS.moduli}

\subsection{Plane sextics and $K3$-surfaces}\label{S.K3}
A \emph{rigid isotopy} of plane projective algebraic curves is
an equisingular deformation or,
equivalently, an isotopy in the class of algebraic curves.
Since, in this paper, we deal with
simple singularities only, the choice of a category (topological,
smooth, piecewise linear) for this definition is irrelevant.
Indeed, recall that one of the fifteen definitions of simple
singularities, see~\cite{Durfee},
is that they are $0$-modal, \ie,
their differential type is determined by their topological type.

Let $C\subset\Cp2$ be a reduced sextic with simple singular
points. Consider the following diagram:
$$
\CD
X   @<<<\BX\\
@Vp VV    @V\Bp VV\\
\Cp2@<\pi<<\BY\rlap,
\endCD
$$
where $X$ is the double covering of~$\Cp2$ branched at~$C$, $\BX$
is the minimal resolution of singularities of~$X$, and $\BY$ is
the minimal embedded resolution of singularities of~$C$ such that
all odd order components of the divisorial pull-back $\pi^*C$
of~$C$ are nonsingular and disjoint. It is well known that $X$ is
a singular $K3$-surface and that $\BX$ is a double covering
of~$\BY$ ramified at
the union of the odd order components of $\pi^*C$.

Let $L\X=H_2(\BX)$; it is a lattice
isomorphic to $2\bE_8\oplus3\bU$. (In what follows we
identify the homology and cohomology of~$\BX$ \via\ the Poincar\'e
duality isomorphism.)
Introduce the following vectors and sublattices:
\Dashes\roster
\dash
$\Gs\X\subset L\X$, the set of the classes of the
exceptional divisors appearing in the blow-up $\BX\to X$;

\dash
$\SS\X\subset L\X$, the sublattice generated by~$\Gs\X$;

\dash
$h\X\in L\X$, the pull-back of the hyperplane section class
$[\Cp1]\in H_2(\Cp2)$;

\dash
$S\X=\SS\X\oplus\<h\X>\subset L\X$;

\dash
$\tSS\X\subset\tS\X\subset L\X$, the primitive hulls
of~$\SS\X$ and~$S\X$, respectively;

\dash
$\Go\X\subset L\X\otimes\R$, the oriented $2$-subspace
spanned by the real and imaginary parts of the class of a
holomorphic $2$-form on~$\BX$ (the `period' of~$\BX$).

\endroster
Clearly, the isomorphism class of the collection $(L\X,h\X,\Gs\X)$
is both a deformation invariant of curve~$C$ and a topological
invariant of pair $(\BY,\pi^*C)$; it is called the
\emph{homological type} of~$C$. By an isomorphism between
two collections $(L',h',\Gs')$ and $(L'',h'',\Gs'')$ we mean an
isometry $L'\to L''$ taking~$h'$ and~$\Gs'$ onto~$h''$
and~$\Gs''$, respectively.

Recall that $\Go\X$ is a positive definite subspace and that the
Picard group $\Pic\BX$ can be
identified with the lattice
$\Go\X^\perp\cap L\X$. In
particular, $\Go\X\in\tS\X^\perp\otimes\R$. Recall also that the
K\"ahler cone~$V^+\X$ of~$\BX$ can be given by
$$
V^+\X=\bigl\{x\in\Go\X^\perp\bigm|
 \text{$x^2>0$ and $x\cdot[E]>0$ for any $(-2)$-curve~$E\subset\BX$}\bigr\}.
$$
The projectivization $\PP(V^+_X)$ is one of the (open) fundamental
polyhedra of the group of motions of the hyperbolic space
$\PP(\{x\in\Go_X^\perp\,|\,x^2>0\})$ generated by the reflections
defined by the roots of $\Pic\BX$.
The walls
bounding~$V^+\X$ are precisely those defined by the classes of the
irreducible $(-2)$-curves in~$\BX$, and the integral classes in the closure
$\barV^+\X$ are the numerically effective divisors on~$\BX$.

In particular, $\Gs\X$ is a `standard' basis of
the root system~$\SS\X$,
so that the cone
$$
W\X=\bigl\{x\in\SS\X\otimes\R\bigm|
 \text{$x\cdot r>0$ for each $r\in\Gs\X$}\bigr\}
$$
is a Weyl chamber of~$\SS\X$. Clearly, $W\X$ and $\Gs\X$
determine each other.

\Remark
Instead of the oriented real subspace~$\Go\X$ one often considers
the Hodge subspace $H^{2,0}(X)\subset L\X\otimes\C$ or,
equivalently, the class $\Go_\C\in L\X\otimes\C$ of a holomorphic
$2$-form on~$\BX$, the latter being defined up to a nonzero factor
and satisfying the conditions $\Go_\C^2=\bar\Go_\C^2=0$,
$\Go_\C\cdot\bar\Go_\C>0$. Then $\Go\X$ is the real part of the
space $H^{2,0}\oplus\,\overline{\!H^{2,0}\!}\,$, or, equivalently,
$\Go\X$ is spanned by the real part~$\Re\Go_\C$ and imaginary
part~$\Im\Go_\C$.
Conversely, $\Go_\C$ can be recovered as $x+iy$, where $x,y$ is
any positively oriented orthonormal basis for~$\Go\X$.
\endRemark

\subsection{Homological types}\label{S.homology}
The set of simple singularities of a plane sextic $C\subset\Cp2$
can be viewed as
a root system~$\SS$ with a distinguished `standard' basis~$\Gs$,
or, equivalently, a distinguished Weyl chamber
$$
W=\bigl\{x\in\SS\otimes\R\bigm|
 \text{$x\cdot r>0$ for each $r\in\Gs$}\bigr\}.
$$
Similar to Section~\ref{S.K3},
let $S=\SS\oplus\<h>$, $h^2=2$.
One has
$\CS=\discr\SS\oplus\<\frac12>$.
An isometry of~$S$ is called \emph{admissible} if it preserves
both~$h$ and~$\Gs$ (as a set).
The group $\Oh S\subset\OAut S$ of admissible isometries
is the group of symmetries of the distinguished Weyl chamber~$W$.
Hence, its image
$\Aut_h\CS\subset\Aut\CS$ coincides with the subgroup $\Aut_\SS\discr\SS$,
see Proposition~\ref{root.auto}. In
particular, $\Aut_h\CS$ is independent of the choice of~$\Gs$.

\definition\label{def.configuration}
Let $\SS$ and~$h$ be as above.
A \emph{configuration} is a finite index
extension $\tS\supset S=\SS\oplus\<h>$ satisfying the following
conditions:
\roster
\item\local1
$\gR(\tSS)=\SS$, where
\smash{$\tSS=h^\perp_{\tS}$}
is the primitive hull
of~$\SS$ in~$\tS$
and $\gR(\tSS)\subset\tSS$ is the sublattice generated by the roots
of \smash{$\tSS$}, see~\ref{S.notation};
\item\local2
there is no root $r\in\SS$ such that $\frac12(r+h)\in\tS$.
\endroster
An isometry of a configuration~$\tS$ is called
\emph{admissible} if it preserves~$S$ and induces an admissible
isometry of~$S$.
\enddefinition

The group of admissible isometries
of~$\tS$ and its image in $\Aut\tCS$ are denoted by $\Oh\tS$
and $\Aut_h\tCS$, respectively. Since $\SS=\gR(\tSS)$ is a
characteristic sublattice of
\smash{$\tSS=h^\perp_{\tS}$},
any
isometry of~$\tS$ preserving~$h$ preserves~$\SS$. Hence, one
has $\Aut_h\tCS=\bigl\{s\in\Aut_h\CS\mid s(\CK)\subset\CK\bigr\}$,
where $\CK$ is the kernel of the extension $\tS\supset S$.

\definition\label{def.AHT}
An \emph{abstract homological type} (extending a fixed set
of simple singularities $(\SS,\Gs)$) is an extension of the
orthogonal sum
$S=\SS\oplus\<h>$, $h^2=2$, to a lattice~$L$ isomorphic to
$2\bE_8\oplus3\bU$ so that the primitive hull~$\tS$ of~$S$ in~$L$
is a configuration.
An \emph{isomorphism} between two abstract homological types
$L_i\supset S_i\supset\Gs_i\cup\{h_i\}$, $i=1,2$, is an
$\Oh S$-isomorphism of extensions, see
Section~\ref{S.extension}; in other words, it is an isometry
$L_1\to L_2$ taking~$h_1$ to~$h_2$ and $\Gs_1$ to~$\Gs_2$ (as a set).
\enddefinition

An abstract homological type is uniquely determined by
the collection $(L,h,\Gs)$; then $\SS$ is the sublattice spanned
by~$\Gs$, and $S=\SS\oplus\<h>$.
The lattices~$\tSS$ and~$\tS$ are defined as the respective
primitive hulls.

\definition\label{def.orientation}
An \emph{orientation} of an abstract homological type
$\CH=(L,h,\Gs)$ is a choice of one of the two orientations of
positive definite $2$-subspaces of the space $\tS^\perp\otimes\R$.
(Recall that $\Gs_+\tS^\perp=2$ and, hence, all
positive definite $2$-subspaces of $\tS^\perp\otimes\R$ can be
oriented in a coherent way.)
The type~$\CH$ is called
\emph{symmetric} if $(\CH,\theta)$ is isomorphic to
$(\CH,-\theta)$ (for some orientation~$\theta$ of~$\CH$).
In other words, $\CH$ is symmetric if it has an automorphism
whose
restriction to~\topsmash{$\tS^\perp$} is $+$-disorienting.
\enddefinition

\subsection{Marked sextics}\label{S.marked}
Let $(\SS,\Gs)$ and $S=\SS\oplus\<h>$ be as in
Section~\ref{S.homology}, and fix an extension $L\supset S$ with
$L\cong2\bE_8\oplus3\bU$. A \emph{marking} (more precisely,
\emph{$(L,h,\Gs)$-marking}) of a singular plane sextic
$C\subset\Cp2$ is an isometry $\Gf\:L\X\to L$ taking~$h\X$
and~$\Gs\X$ onto~$h$ and~$\Gs$, respectively (see
Section~\ref{S.K3} for the notation). A \emph{marked} sextic is a
sextic supplied with a distinguished marking.

The following statement, based on the surjectivity of the period
map for $K3$-surfaces, is essentially contained in T.~Urabe~\cite{Urabe}.

\proposition\label{existence}
Let $(L,h,\Gs)$ be a collection as above, and let~$\Go$ be
an oriented positive definite $2$-subspace in
$S^\perp_L\otimes\R$.
Then there exists a singular plane sextic $C\subset\Cp2$ and an
$(L,h,\Gs)$-marking $\Gf\:L\X\to L$ taking~$\Go\X$ to~$\Go$
if and only if the following conditions are satisfied\rom:
\roster
\item\local1
$(L,h,\Gs)$ is an abstract homological type\rom;
\item\local2
every root $r\in L$ orthogonal to both~$h$ and~$\Go$ belongs
to~$\SS$.
\endroster
\endproposition

We precede the proof with a lemma. Denote by~$\Gamma$ the group
generated by the reflections defined by the roots of the lattice
$\Go^\perp\cap L$.

\lemma\label{V.unique}
Let $(L,h,\Gs)$ and $\Go$ be as in
Proposition~\ref{existence}, and assume that
conditions~\iref{existence}1,~\ditto2 are satisfied.
Then there is a unique open
convex cone $V^+=V^+(\Go)\subset\Go^\perp$ such that
\Dashes\roster
\dash
the projectivization $\PP(V^+_X)$ is one of the fundamental
polyhedra of the action of~$\Gamma$ on the hyperbolic space
$\PP(\{x\in\Go^\perp\,|\,x^2>0\})$\rom;
\dash
the closure $\barV^+$ contains~$h$\rom;
\dash
the intersection $V^+\cap(\SS\otimes\R)$ is the Weyl chamber~$W$
defined by~$\Gs$.
\endroster
\endlemma

\proof
Condition~\iref{existence}2
implies that $W$ extends to a Weyl
chamber~$W'$
in the negative definite space
$h^\perp_{\smash{\Go^\perp}}$; it is characterized by the
requirement that $W'\cdot r>0$ for each $r\in\Gs$. Then, Vinberg's
algorithm~\cite{Vinberg} applied to~$h$ extends~$\PP(W')$ to a
unique fundamental polyhedron~$P$ of~$\Gamma$ whose
closure~$\overline{P}$ contains the class $h/\R^*$. The
connected component of the cone $\{x\in\Go^\perp\,|\,x/\R^*\in P\}$
containing~$W$
is the desired
cone~$V^+$.
\endproof

\proof[Proof of Proposition~\ref{existence}]
In the presence of~\itemref{existence}2, condition~\ditto1 is
equivalent to the requirement that
\roster
\item[3]\local3
there is no element $u\in\Go^\perp\cap L$ with $u^2=0$ and
$u\cdot h=1$.
\endroster
In this form, it is obvious that the conditions are necessary:
\loccit3 is necessary for the linear system~$h$ to define a
degree~$2$ map $\BX\to\Cp2$, see~\cite{Urabe}, and \ditto2 means
that the curves contracted by this map are exactly those defined
by the elements of~$\Gs$, \ie, the sextic does have the prescribed
set of singularities.

Prove the sufficiency. Due to the surjectivity of the period map,
there is a $K3$-surface~$\BX$ and an isomorphism $\Gf\:H_2(\BX)\to
L$ taking~$\Go\X$ to~$\Go$. The image $\Gf(V^+\X)$ of the the
K\"ahler cone~$V^+\X$ of~$\BX$ is a fundamental domain of the
action of~$\Gamma$ on one of the two halves of the positive cone
$\{x\in\Go^\perp\,|\,x^2>0\}$. Hence, composing~$\Gf$ with an
element of~$\Gamma$ and, if necessary, multiplication by~$-1$, one
can assume that $\Gf$ takes~$V^+\X$ to the cone~$V^+(\Go)$ given
by Lemma~\ref{V.unique}. Then the pull-back $h\X=\Gf^{-1}(h)$
belongs to the closure~$\barV\X$; hence, it is numerically
effective and, due to condition~\loccit3, it defines a degree~$2$ map
$p\:\BX\to\Cp2$, see~\cite{Urabe}. The elements of the pull-back
$\Gs\X=\Gf^{-1}(\Gs)$ define (some of) the walls of the K\"ahler
cone and, hence, are realized by irreducible $(-2)$-curves
in~$\BX$; due to condition~\ditto2, they are all the $(-2)$-curves
contracted by~$p$. Thus, $\Gf$ is the desired marking.
\endproof

\subsection{Moduli}\label{S.moduli}
In view of Proposition~\ref{existence}, when speaking about
$(L,h,\Gs)$-marked sextics, one can assume that $\CH=(L,h,\Gs)$ is
an abstract homological type. Since the period~$\Go\X$ changes
continuously within a family, the orientation of the image
$\Gf(\Go\X)$ is an additional discrete invariant of deformations
in the class of marked plane sextics.

\theorem\label{2.components}
For each abstract homological type $\CH=(L,h,\Gs)$ there
are exactly two rigid isotopy classes of $\CH$-marked plane sextics.
They differ by the orientation of the positive definite
$2$-subspace $\Gf(\Go\X)\subset\tS^\perp\otimes\R$.
\endtheorem

\proof
The existence of at least two rigid isotopy classes that differ by
the orientation of $\Gf(\Go\X)$ is given by
Proposition~\ref{existence}. Thus,
it suffices
to show that any two $\CH$-marked $K3$-surfaces
$(\BX_0,\Gf_0)$, $(\BX_1,\Gf_1)$ satisfying~\iref{existence}2 and such
that the images $\Gf_t(\Go\X)$, $t=0,1$, have coherent
orientations can be connected by a family $(\BX_t,\Gf_t)$,
$t\in[0,1]$ of $\CH$-marked $K3$-surfaces still
satisfying~\iref{existence}2. Then the linear systems
$h_t=\Gf_t^{-1}(h)$ would define a family of degree~$2$ maps
$\BX_t\to\Cp2$
and, since \iref{existence}2 holds for each~$t$, the resulting
family $C_t\in\Cp2$ of the branch curves would be equisingular.

Consider the space $\GOK$ of pairs $(\Go,\Gr)$, where
$\Go\subset L\otimes\R$ is an oriented positive definite $2$-subspace
and $\Gr\in L\otimes\R$ is a positive vector ($\Gr^2>0$) orthogonal
to~$\Go$. Let
$$
\GOK_0=\GOK\sminus\bigcup_{r\in L,\,r^2=-2}
 \bigl\{(\Go,\Gr)\in\GOK\bigm|\Go\cdot r=\Gr\cdot r=0\bigr\}.
$$
According to A\.~Beauville~\cite{Beauville}, $\GOK$ is a fine
period space
of marked quasi-polarized $K3$-surfaces, a
\emph{quasi-polarization} being a class
of a K\"ahler metric.

Let $\GO(\CH)\cong\O(2,d)/\SO(2)\times\O(d)$ be the space of
oriented positive definite $2$-subspaces $\Go\subset
S^\perp\otimes\R$ (here $d=19-\rank S$), and let
$\GO_0(\CH)\subset\GO(\CH)$ be the set of subspaces~$\Go$
satisfying~\iref{existence}2. Since $\CH$ is an abstract
homological type, $\GO_0(\CH)$ is obtained from $\GO(\CH)$ by
removing a countable number of codimension~$2$ subspaces
$H_r=\{\Go\mid\Go\cdot r=0\}$, $r\in h^\perp\sminus\SS$, $r^2=-2$.
Condition~\iref{def.configuration}1 implies that none of~$H_r$
coincides with $\GO(\CH)$ and, hence, $\GO_0(\CH)$ is nonempty.
Since $\GO(\CH)$ has two connected components, so does
$\GO_0(\CH)$. The components differ by the orientation of the
subspaces.

Now, let $\GOK_0(\CH)\subset\GOK_0$ be the subspace
$\{(\Go,\Gr)\mid\Go\in\GO_0(\CH),\ \Gr\in V^+(\Go)\}$, where
$V^+(\Go)$ is the
cone given by Lemma~\ref{V.unique}. In
view of Proposition~\ref{existence} and Lemma~\ref{V.unique},
Beauville's result cited above implies that $\GOK_0(\CH)$ is a fine
period space of $\CH$-marked quasi-polarized plane sextics. On the
other hand, the natural projection $\GOK_0(\CH)\to\GO_0(\CH)$,
$(\Go,\Gr)\mapsto\Go$, has contractible fibers (the cones
$V^+(\Go)$) and, outside a countable union of codimension~$2$
subsets $H_r$, $r\in L\sminus\SS$, $r^2=-2$, it is a locally
trivial fibration. Hence, the period space $\GOK_0(\CH)$ has two
connected components, and the statement follows.
\endproof

\theorem\label{enumeration}
The map sending a plane sextic $C\subset\Cp2$ to the pair
consisting of its homological type $(L\X,h\X,\Gs\X)$ and the
orientation of the space~$\Go\X$ establishes a one-to-one
correspondence between the set of rigid isotopy classes of sextics
with a given set of singularities $(\SS,\Gs)$ and the set of
isomorphism classes of oriented abstract homological types
extending~$(\SS,\Gs)$.
\endtheorem

\proof
The statement is an immediate consequence of
Theorem~\ref{2.components} and the obvious fact that any two
$\CH$-markings of a given sextic differ by an isometry
of the abstract homological type~$\CH$.
\endproof

\Remark
Marked sextics satisfying conditions~\iref{existence}1 and~\ditto2
can be regarded as ample marked $\tS$-polarized $K3$-surfaces in
the sense of Nikulin \cite{NikK3}. (The ampleness of the
polarization follows from condition~\iref{existence}2.) Their period space
is constructed in Dolgachev~\cite{Dolgachev}, based directly on the global
Torelli theorem and the surjectivity of the period map. It is
shown that the period space has two connected components
interchanged by the complex conjugation.
\endRemark

\subsection{Proof of Theorem~\ref{th.main}}\label{main.proof}
The `only if' part of the statement is obvious. We will prove the
`if' part under the assumption that $C_1$ has at least one singular
point. (Otherwise the two sextics are nonsingular and, hence,
rigidly isotopic.)

The regularity condition implies that $f$ preserves the complex
orientations of both~$\Cp2$ and~$C_1$; in particular, the induced
map $f_*\:H_2(\Cp2)\to H_2(\Cp2)$ takes~$[\Cp1]$ to~$[\Cp1]$.
Furthermore, $f$ lifts to
a diffeomorphism $\BY_1\to\BY_2$ and, hence, to a diffeomorphism
$\Bf\:\BX_1\to\BX_2$ of the corresponding $K3$-surfaces (see
Section~\ref{S.K3} for the notation). The induced homomorphism
$\Bf_*\:L_{X_1}\to L_{X_2}$ takes~$h_{X_1}$ and~$\Gs_{X_1}$
to~$h_{X_2}$ and~$\Gs_{X_2}$, respectively. Hence, for each
marking $\Gf\:L_{X_2}\to L$ of~$C_2$ the composition
$\Gf\circ\Bf_*$ is a marking of~$C_1$. The crucial observation is
the fact that, according to S\.~K\.~Donaldson~\cite{Donaldson},
the map~$\Bf_*$ induced by a \emph{diffeomorphism} of
$K3$-surfaces preserves the orientation of the (positive
definite) $3$-subspace spanned by the period $\Go_{X_1}$ and a
K\"ahler class~$\Gr_{X_1}$. Since K\"ahler classes~$\Gr_{X_1}$
and~$\Gr_{X_2}$ can be chosen arbitrary close to~$h_{X_1}$
and~$h_{X_2}$, respectively (recall that
$h_{X_1}$
and~$h_{X_2}$ belong to the closures of the respective K\"ahler
cones), the latter assertion means that the orientations of
$\Gf(\Go_{X_2})$ and $\Gf\circ\Bf_*(\Go_{X_2})$ agree, and
Theorem~\ref{2.components} implies that $C_1$ and~$C_2$ are
rigidly isotopic.
\qed

\section{Geometry of plain sextics}\label{SS.sextics}

In this section, we discuss the relation between the geometry of a
plane sextic and its homological type. We start with introducing
several versions of the notion of Zariski pair
(Section~\ref{S.Zariski.def}) and outlining Yang's algorithm
recovering the combinatorial data of a curve from its
configuration (Section~\ref{S.Yang}).
Sections~\ref{S.reducible}
and~\ref{S.Zariski.geom} give a simple characterization of,
respectively, reducible and abundant sextics.

\subsection{Zariski pairs}\label{S.Zariski.def}
Two reduced curves $C_1,C_2\subset\Cp2$ are said to have the same
\emph{combinatorial data} if there exist irreducible decompositions
$C_i=C_{i,1}+\ldots+C_{i,k_i}$, $i=1,2$, such that:
\roster
\item
$k_1=k_2$ and $\deg C_{1,j}=\deg C_{2,j}$ for all
$j=1,\ldots,k_1$;
\item
there is a one-to-one correspondence between the singular points
of~$C_1$ and those of~$C_2$ preserving the topological types of
the points;
\item
two singular points $P_i\in C_i$, $i=1,2$,
corresponding to each other are related by a local
homeomorphism such that if a branch at~$P_1$ belongs to a
component~$B_{1,j}$ then its image belongs to~$B_{2,j}$.
\endroster
For an irreducible curve~$C$, its combinatorial data are determined
by the degree $\deg C$ and the set of topological types of the
singularities of~$C$.

One of the principal question in topology of plane curves is the
extent to which the combinatorial data of a curve determine its
global behavior. In order to formalize this question,
Artal~\cite{Artal.def} suggested the notion of Zariski pair.

\definition\label{def.Zariski}
Two reduced curves $C_1,C_2\subset\Cp2$ are said to form a
\emph{Zariski pair} if
\roster
\item\local1
$C_1$ and~$C_2$ have the same combinatorial data, and
\item\local2
the pairs $(\Cp1,C_1)$ and $(\Cp2,C_2)$ are not homeomorphic.
\endroster
\enddefinition

\Remark\label{rem.Zariski}
Cited above is the more suitable
definition used in the subsequent papers. The original definition
suggested in~\cite{Artal.def} requires, instead of
\iref{def.Zariski}1, that the pairs $(T_1,C_1)$ and $(T_2,C_2)$
should be diffeomorphic, where $T_i$ is a regular neighborhood
of~$C_i$, $i=1,2$. If the singularities involved are simple, the
two definitions are equivalent
as, on the one hand, simple singularities are $0$-modal and, on
the other hand, simple curve singularities are distinguished by
their links (which is a straightforward consequence of their
classification).
\endRemark

Condition~\itemref{def.Zariski}2 in Definition~\ref{def.Zariski}
varies from paper to paper: one can replace it with the negation
of any reasonable `global' equivalence relation. For example,
instead of~\ditto2 it is sometimes required that the complements
$\Cp2\sminus C_1$ and $\Cp2\sminus C_2$ should not be
homeomorphic. Relevant for the present paper are the following
notions:
\Dashes\roster
\dash
\emph{regular Zariski pair}, with~\iref{def.Zariski}2 replaced by
`the pairs $(\Cp1,C_1)$ and $(\Cp2,C_2)$ are not regularly
diffeomorphic in the sense of Theorem~\ref{th.main};'
\dash
\emph{classical Zariski pair}, with~\iref{def.Zariski}2 replaced
by `the Alexander polynomials $\Delta_{C_1}(t)$ and
$\Delta_{C_2}(t)$ differ,' see Section~\ref{S.Zariski.geom} for
details.
\endroster
Theorems~\ref{th.main} and~\ref{enumeration} state that, in order
to construct examples of regular Zariski pairs, it suffices to
find curves with the same combinatorial data but not isomorphic
oriented homological types. The notion of classical Zariski pair
is of a historical interest, as
it was the Alexander polynomial that was used to distinguish the
curves in the first examples.
In Section~\ref{S.Zariski} below we enumerate the deformation
families of unnodal curves
whose Alexander polynomial is not determined by the combinatorial
data.

\subsection{Configurations and combinatorial data}\label{S.Yang}
Let $C_1,C_2\subset\Cp2$ be a pair of reduced plane sextics with
simple singularities.
Consider the corresponding oriented homological types
$(\CH_i,\theta_i)=(L_i,h_i,\Gs_i,\theta_i)$, $i=1,2$, and related
lattices $\SS_i$, $\tS_i$, \etc., see Section~\ref{S.K3}. (To
simplify the notation we use index~$i$ instead of~$X_i$.) Recall
that the finite index extension $\tS\supset\SS\oplus\<h>$ is
called a configuration, see Definition~\ref{def.configuration}.

\midinsert
$$
\def\bx#1{\vcenter{\hsize.4\hsize\parindent0pt\normalbaselines\raggedright#1}}
\CD
\SS_1\cong\SS_2@:<===>
 \bx{$C_1$ and $C_2$ have the same set of singularities}\\
@+<|\double\up|>@+<|\double\up|>\\
(\tS_1,h_1,\Gs_1)\cong(\tS_2,h_2,\Gs_2)@:====>
 \bx{$C_1$ and $C_2$ have the same combinatorial data;\endgraf
 $\Delta_{C_1}(t)=\Delta_{C_2}(t)$}\\
@+<|\double\up|>@+<|\double\up|>\\
\CH_1\cong\CH_2@:<===>
 \bx{$C_1$ is rigidly isotopic to either $C_2$ or its conjugate
 $\bar C_2$}\\
@+<|\double\up|>@+<|\double\up|>\\
(\CH_1,\theta_1)\cong(\CH_2,\theta_2)@:<===>
 \bx{$C_1$ is rigidly isotopic to $C_2$;\endgraf
 $(\Cp2;C_1)$ and $(\Cp2;C_2)$ are regularly diffeomorphic}
\endCD
$$
\figure[Diagram]
\endfigure\label{arithmetics=>geometry}
\endinsert

Diagram~\ref{arithmetics=>geometry} represents various relations
between
the geometric properties of a curve and
arithmetic properties of its homological type.
The equivalence in the first line of the diagram is obvious; the
equivalences in the last two lines are the statement of
Theorem~\ref{enumeration} and the fact that the two components of
the period space are interchanged by the complex conjugation.

Informally, the implication in the second line of the diagram
states that the configuration~$\tS_X$ encodes the existence of
various auxiliary curves passing in a prescribed way through the
singular points of \emph{each} curve deformation equivalent
to~$C$. In short, this assertion follows from the fact that
$\tS_X$ is the Picard group of a generic curve of the deformation
family. The precise algorithm recovering the combinatorial data of
a sextic~$C$ from its configuration~$\tS_X$ is outlined at the end
of this section. The relation between the configuration and the
Alexander polynomial in the case of irreducible curves is
discussed in Section~\ref{S.Zariski.geom}.

Note that the implication in the second line of
Diagram~\ref{arithmetics=>geometry} is not invertible. There are
pairs of curves with the same combinatorial data and/or Alexander
polynomial but not isomorphic configurations, see, \eg,
Theorem~\ref{th.Zariski} and Proposition~\ref{2A9+A1}.

\subsubsection*{Yang's algorithm}
Assume that a sextic~$C$ splits into irreducible
components $C_1,\ldots,C_k$. Consider the fundamental classes
$[C],[C_i]\in H^2(\BX)$ in the homology of the covering
$K3$-surface. They realize certain elements $c,c_i$ of the group
$\tS_X\subset S_X^*$, so that $c=\sum_ic_i$. The classes $c$,
$c_i$ are recovered from the combinatorial data of~$C$: one has
$c\cdot h_X=6$, $c_i\cdot h_X=\deg C_i$, and the intersections
$c\cdot\Gs_j$, $c_i\cdot\Gs_j$ with the classes of the exceptional
divisors are determined by the incidence of the curves in the
minimal resolution of singularities. In fact, to each local
branch~$b$ at a simple singular point~$P$ one can assign an element
$\che(b)$ of the group $(\SS_P)^*$ dual to the lattice~$\SS_P$
spanned by the exceptional divisors at~$P$. Then, for a
component~$C_i$ of~$C$, one has
$$
c_i=\frac12(\deg C_i)h_X+\sum_{b\in C_i}\che(b).\eqtag\label{eq.splitting}
$$
Explicit expressions for the elements~$\che(b)$ are found
in~\cite{Yang}.
Next lemma is an immediate consequence of these formulas.

\lemma\label{order2}
If a simple singular point~$P$ has more than one branch $b_1,\ldots,b_k$,
then each residue $\che(b_i)\bmod\SS_P\in\discr\SS_P$ is an
element of order~$2$\rom;
these residues are subject to
the only relation $\sum_{i=1}^k\che(b_i)=0\bmod\SS_P$.
\endlemma

Now, one can ignore the geometric setting and consider a
\emph{virtual decomposition}, \ie, a decomposition $c=\sum_ic_i$
determined by a hypothetical set of combinatorial data of~$C$.
Certainly, \emph{a priori} one can only assert that $c_i\in
S_X^*$. The set of all virtual decompositions of~$c$ is partially
ordered by degeneration, so that $c=c$ is the minimal element. The
following statement is based on the Riemann-Roch theorem for
$K3$-surfaces.

\theorem[Theorem \rm(see Yang~\cite{Yang})]\label{th.Yang}
The actual combinatorial data of an irreducible plane sextic~$C$
with simple singularities is the one corresponding to the only
maximal element in the set of the virtual decompositions
$c=\sum_ic_i$ with all $c_i\in\tS_X$.
\endtheorem

\subsection{Reducible curves}\label{S.reducible}
In this section, $C$ is a reduced plane curve of arbitrary degree
$d=4m+2$. We still assume that all singularities of~$C$ are
simple. As in the case of sextics, consider the double
covering~$X$ branched over~$C$ and its minimal resolution~$\BX$.
Certainly, $\BX$ is not a $K3$-surface; however, since $\BX$ is
diffeomorphic to the double covering branched over a nonsingular
curve, one still has $\pi_1(\BX)=0$ and $L_X=H_2(\BX)$ is an even
lattice.

All lattices $\SS_X\subset S_X=\SS_X\oplus\<h_X>\subset\tS_X\subset L_X$
introduced in Section~\ref{S.K3} for sextics still make sense in
the general case.

\theorem\label{th.reducible}
Let~$C$ be a reduced plane curve of degree $4m+2$ and with simple
singularities only. Then $C$ is reducible if and only if the
kernel~$\CK$ of the extension $\tS_X\supset S_X$ has elements of
order~$2$.
\endtheorem

\proof
If $C_i$ is a proper component of~$C$,
the residue $c_i\bmod S_X\in\CK$ given by~\eqref{eq.splitting} is an
element of order~$2$. It is nontrivial since $C_i$ must pass
through a singular point of~$C$ that is not entirely contained
in~$C_i$. (Clearly, the calculation of~$c_i$,
including Lemma~\ref{order2}, still applies to the general case of
curves of degree $4m+2$.)

Now, assume that $C$ is irreducible.
Denote by~$\BC\subset\BX$ the set theoretical pull-back of~$C$; it
is the union of the exceptional divisors and the proper pull-back,
which can be identified with~$C$ itself. Consider the fundamental
group $\pi=\pi_1(\Cp2\sminus C)$ and the homomorphism
$\kappa\:\pi\to\Z/2\Z$ defining the double covering
$\BX\sminus\BC\to\Cp2\sminus C$. One has
$\Ker\kappa=\pi_1(\BX\sminus\BC)$.

The abelinization $\pi/[\pi,\pi]=H_1(\Cp2\sminus C)$ is the cyclic
group $\Z/(4m+2)\Z$; its $2$-primary part is~$\Z/2\Z$, and
from the exact sequence
$$
\{1\}@>>>\Ker\kappa/(\Ker\kappa)^2@>>>\pi/(\Ker\kappa)^2@>>>\Z/2\Z@>>>\{1\}
$$
and properties of $2$-groups one concludes that
the group $\Ker\kappa/(\Ker\kappa)^2=H_1(\BX\sminus\BC;\F_2)$ is
trivial. Then, from the Poincar\'e duality and the fact that
$H^3(\BX;\F_2)=0$ it follows that
the inclusion homomorphism $H^2(\BX;\F_2)\to H^2(\BC;\F_2)$ is onto and
its dual $H_2(\BC;\F_2)\to H_2(\BX;\F_2)$ is a monomorphism.
On the other hand, since $C$ is
irreducible and $[C]=(2m+1)h_X\bmod\SS_X$,
the inclusion induces an isomorphism
$H_2(\BC;\F_2)=S_X\otimes\F_2$. Thus, the $(\bmod2)$-reduction
$S_X\otimes\F_2\to\tS_X\otimes\F_2$ is a monomorphism. This fact
implies that
$\CK$ has no elements of order~$2$.
\endproof

\Remark
In the case of sextics, Theorem~\ref{th.reducible} can as well be
deduced from Theorem~\ref{th.Yang}. However, this would require a
thorough analysis of a number of exceptional cases and eliminating
them using conditions~\iref{def.configuration}1 and~\ditto2 and an extended
version of Proposition~\ref{extension.roots}.
\endRemark

\subsection{Classical Zariski pairs}\label{S.Zariski.geom}
The \emph{Alexander polynomial} of a degree~$m$
plane curve $C\subset\Cp2$
can be defined as the characteristic polynomial of the
deck translation action on $H_1(\BX_m;\C)$, where $\BX_m$ is the
desingularization of the $m$-fold cyclic covering of~$\Cp2$
ramified at~$C$
(see A\.~Libgober~\cite{Libgober.first}--\cite{Libgober.last}
for details).
By a \emph{classical Zariski pair} we mean a pair of curves
that have the same
combinatorial data and differ by their Alexander polynomial. (The truly
classical Zariski pair, due to Zariski himself~\cite{Zariski}, is
a pair of irreducible sextics with six cusps each, one of them
having all cusps on a conic, and the other one not.)

The Alexander polynomials $\Delta_C(t)$ of all irreducible
sextics~$C$ are found in~\cite{LOMI} (see also~\cite{poly}),
where it is shown that $\Delta_C(t)=(t^2-t+1)^d$ and the
exponent~$d$ is determined by the set of singularities of~$C$
unless the latter has the form
$$
\SS=e\bE_6\oplus\bigoplus_{i=1}^6a_i\bA_{3i-1}\oplus n\bA_1,\quad
2e+\sum ia_i=6.\eqtag\label{Zariski.set}
$$
If the set of singularities is as in~\eqref{Zariski.set}, then $d$
may \emph{a priori} take values~$0$ or~$1$; in the latter case the
curve is called \emph{abundant}. The following statement is proved
in~\cite{poly}.

\theorem\label{th.abundant.geom}
For an irreducible plane sextic~$C$ with a set of
singularities~$\SS$ as in~\eqref{Zariski.set}, the following three
conditions are equivalent\rom:
\roster
\item\local1
$C$ is abundant\rom;
\item\local2
$C$ is \emph{tame}, \ie, it is given by an equation of the form
$f_2^3+f_3^2=0$, where $f_2$ and~$f_3$ are some
polynomials of degree~$2$ and~$3$, respectively\rom;
\item\local3
there is a
conic~$Q$ whose local intersection index with~$C$ at
each singular point of~$C$ of type~$\bA_{3i-1}$
\rom(respectively,~$\bE_6$\rom) is~$2i$
\rom(respectively,~$4$\rom).
\endroster
\endtheorem

Observe that the discriminant group of each lattice~$\bA_{3i-1}$
or~$\bE_6$ has a unique subgroup isomorphic to~$\Z/3\Z$. Its
nontrivial elements are the residues of~$\cheb{1,2}$,
where $\cheb1$ is the element given in some standard basis
$e_1,e_2,\ldots$ by
$$
\gather
\cheb1=\frac13(2e_1+4e_2+\ldots+2ie_{i}+(2i-1)e_{i+1}+\ldots+e_{2i-1})
 \in(\bA_{2i-1})^*,\\
\cheb1=\frac13(4e_1+5e_2+6e_3+4e_4+2e_5+3e_6)\in(\bE_6)^*,
\endgather
$$
and $\cheb2$ is obtained from~$\cheb1$ by the only nontrivial
symmetry of the
Dynkin graph. (In the case~$\bE_6$, the basis elements are
numbered so that $e_6$ is attached to the short edge of the
graph.) The following theorem characterizes abundant curves in
terms of configurations.

\theorem\label{th.abundant.config}
Let~$C$ be a plane sextic with a set of singularities~$\SS$ as
in~\eqref{Zariski.set}. Then a reduced conic~$Q$ as
in~\iref{th.abundant.geom}3 exists if and only if the kernel~$\CK$
of the extension $\tS_X\supset S_X$ has $3$-torsion. If this is
the case, the $3$-primary part of~$\CK$ is a cyclic group of
order~$3$ generated by a residue of the form
$\sum\cheb{1,2}_i\bmod S_X$, where $\cheb{1,2}_i$ are the elements
defined above and the sum contains exactly one element for each
singular point of~$C$ other than~$\bA_1$.
\endtheorem

\proof
Assume that a conic~$Q$ exists. Resolving the singularities, one
can see that the proper pull-back of~$Q$ in~$\BY$ does not
intersect the branch locus; hence, $Q$ lifts to a pair of rational
curves (possibly, reducible) in~$\BX$. In the homology of~$\BX$,
each of the lifts realizes a class of the form
$q=h+\sum\cheb{1,2}_i$, the summation involving exactly one
element for each singular point other than~$\bA_1$. Hence,
$q\in\tS_X$ and the residue $q\bmod S_X$ is a $3$-torsion element
of~$\CK$. Conversely, if a class~$q$ as above belongs to~$\tS_X$,
the Riemann-Roch theorem implies that $q$ is realized by a
rational curve. Its projection to~$\Cp2$ is a conic~$Q$ as
in~\iref{th.abundant.geom}3.

Show that any element $q\in\CK$ of order~$3$ must be as in the
statement (and hence give rise to a conic~$Q$).
Clearly, $q$ is a linear combination of the residues
$\cheb{\cdot}_i\bmod S_X$. One has $(\cheb{\cdot}_i)^2=-2i/3\bmod2\Z$
for $\cheb{\cdot}_i\in(\bA_{3i-1})^*$ and
$(\cheb{\cdot}_i)^2=2/3\bmod2\Z$ for $\cheb{\cdot}_i\in(\bE_6)^*$.
Since $q$ is isotropic, it must either involve all singular points
of~$C$ other than~$\bA_1$ or else
belong to an orthogonal summand of~$S_X$ of the form $3\bA_2$,
$\bA_5\oplus\bA_2$, $\bA_8$, or $\bE_6\oplus\bA_2$. The latter
possibility is ruled out by condition~\iref{def.configuration}1
and Proposition~\ref{extension.roots}.

If $q'\in\CK$ were another element of order~$3$, $q'\ne\pm q$,
then the sum $q+q'$ would be an order~$3$ element not involving
all singular points. Hence, $\CK$ contains at most two (opposite)
elements of order~$3$.

Finally, if $q\in\CK$ is an element of order~$9$, then either $3q$
is an element of order~$3$ not involving all singularities, or $q$
is the sum of two generators of $\discr(2\bA_8)$, or $q$ is twice a
generator of $\discr\bA_{17}$. In the last two cases $q$ cannot be
isotropic.
\endproof

\corollary\label{Zariski.conf}
Each set of singularities~$\SS$ as in~\eqref{Zariski.set}
extends to two isomorphism classes of
configurations $\tS\supset S=\SS\oplus\<h>$
that may correspond to irreducible sextics, one abundant
\rom($\CK=\Z/3\Z$\rom) and one not \rom($\CK=0$\rom).
\endcorollary

\proof
The $2$-primary part of the kernel~$\CK$ is trivial due to
Theorem~\ref{th.reducible}; the $3$-primary part is given by
Theorem~\ref{th.abundant.config}. All extensions with $\CK=\Z/3\Z$
are isomorphic to each other as the two elements $\cheb{1,2}$
corresponding to each singularity are interchangeable by an
admissible automorphism. Finally, $\CS$ cannot have an isotropic
subgroup of prime order other than~$2$ or~$3$. In fact, the only
nontrivial $p$-primary component, $p\ne2$ or~$3$, is
$\CS\otimes\Z_5\cong\<-\frac25>$ in the case
$\SS=\bA_{14}\oplus\bA_2$. It has no isotropic elements.
\endproof

\Remark\label{rem.Zariski.orientation}
Since elements $\cheb{1,2}$ are not symmetric,
Theorem~\ref{th.abundant.config} implies that the singular points
of an abundant curve admit a natural coherent `orientation'
(order of the exceptional divisors in~$\BX$). Geometrically, this
order is selected by a choice of one of the two components of the
pull-back of~$Q$ in~$\BX$. If the order of the exceptional
divisors were fixed, instead of Corollary~\ref{Zariski.conf} one
would have $2^{m-1}$ non-isomorphic abundant configurations, where
$m$ is the number of the singular points other than~$\bA_1$.
\endRemark

\section{Examples}\label{SS.examples}

Theorem~\ref{enumeration} reduces rigid isotopy classification of
plane sextics to the enumeration of oriented abstract homological
types. In this concluding section, we outline the principal steps
of the classification and illustrate them on a few examples:
sextics with few singularities, where Nikulin's theorems apply to
give a unique rigid isotopy class (Section~\ref{S.few.points}),
unnodal classical Zariski pairs (Section~\ref{S.Zariski}), and a
few recent examples of sextics with maximal total Milnor number
(Section~\ref{S.examples}). Section~\ref{S.remarks} contains a few
concluding remarks, speculations, and open problems.

\subsection{Enumerating abstract homological types}\label{S.arithmetics}
Recall that an isomorphism of abstract homological types is
defined as an isometry preserving the distinguished class~$h$ and
distinguished basis~$\Gs$ (as a set). Next proposition states that
$\Gs$ can be ignored.

\proposition
Let $\CH_i=(L_i,h_i,\Gs_i)$, $i=1,2$, be two abstract homological
types, and let $S_i$ be the corresponding sublattices
spanned by~$h_i$ and~$\Gs_i$. Then $\CH_1$, $\CH_2$
are isomorphic if and only if there is an isometry $t\:L_1\to L_2$
taking $S_1$ to~$S_2$ and $h_1$ to~$h_2$.
\endproposition

\proof
The extensions of the restriction~\smash{$t_{\tS_1^\perp}$} to the whole
lattice~$L_1$ depend only on the induced map $\tCS_1\to\tCS_2$. In
view of Proposition~\ref{root.auto}, the image in~$\Aut\tCS_i$ of
the group of admissible isometries of~$\tS_i$ coincides with the
image of the group of isometries preserving~$h_i$.
\endproof


Fix a set of singularities~$\SS$. The classification of oriented
abstract homological types extending~$\SS$ is done in four steps.

\subsubsection**{Step~1\rom: enumerating the configurations~$\tS$ extending~$\SS$}
Due to Proposition~\ref{extension.finite},
a configuration is determined by a choice of an
isotropic subgroup $\CK\subset\CS$.
Note that, given~\iref{def.configuration}1,
condition~\iref{def.configuration}2 should only
be checked for the direct summands of~$\SS$
isomorphic to~$\bA_1$, as for any other root
$r\in\SS$ there is another root $r'\in\SS$
such that $r\cdot r'=1$ and, hence, $r+h$ is primitive in~$\tS$.
We combine this observation and
Corollary~\ref{extension.roots} to the following statement.

\proposition\label{configuration}
Let~$\tS$ be a configuration extending a set of
singularities~$\SS$. Then
each direct summand of~$\SS$ isomorphic to one of the root systems
listed in Corollary~\ref{extension.roots} is primitive
in~$\tS$, and
each sublattice $\bA_1\oplus\<h>$, where $\bA_1$ is a direct
summand of~$\SS$, is primitive in~\topsmash{$\tS$}.
\endproposition


\subsubsection**{Step~2\rom: enumerating the isomorphism classes of $\tS^\perp$}
The orthogonal complement $N=\tS^\perp_L$ has genus
$(2,19-\rank\SS;-\tCS)$.
The existence of a lattice in this genus, whenever it
holds, is given by Theorem~\ref{N.existence}. If
$N$ is indefinite,
one would hope that a
theorem similar to~\ref{N.uniqueness} would
imply uniqueness.
The case of definite lattices ($\rank\SS=19$)
is treated in Section~\ref{S.rank2}. There are examples
(see, \eg, Proposition~\ref{A18+A1} below) when the genus does
contain more than one isomorphism class.

\subsubsection**{Step~3\rom: enumerating the bi-cosets of
$\Aut_h\tCS\times\Aut_N\CN$}
Once the lattice $N=\tS^\perp$ is chosen, one can
fix an anti-isometry $\tCS\to\CN$ and, hence, an
isomorphism $\Aut\CN=\Aut\tCS$; then, the extensions are
classified by the quotient set
$\Aut_h\tCS\backslash\!\Aut\tCS/\!\Aut_N\CN$. Important special cases
are those with $\Aut_h\tCS=\Aut\tCS$ (\cf.
Section~\ref{S.notation}) or $\Aut_N\CN=\Aut\CN$ (this would
normally be given by Theorem~\ref{N.strong.uniqueness}). If none
of the above applies, the isometries are to be described
manually.
There are
examples (see Propositions~\ref{A19}, \ref{2A9+A1},
and~\ref{E6+A7+A3+A2+A1} below) when the quotient consists of more
than one coset, thus giving rise to more than one abstract
homological type.

\subsubsection**{Step~4\rom: detecting whether the abstract
homological types are symmetric}
An abstract homological type is
symmetric if and only if $\tS^\perp$ has a $+$-disorienting
isometry~$t$ whose image in $\Aut\discr\tS^\perp=\Aut\tCS$
belongs to the
product of the subgroup $\Oh\tS$ and the image of $\OPlus{\tS^\perp}$.
Asymmetric abstract homological types do exist, see
Proposition~\ref{A18+A1}.
Below is a sufficient condition for an
abstract homological type to be symmetric.

\proposition\label{symmetric.-2}
Let $\CH=(L,h,\Gs)$ be an abstract homological type. If the
lattice $\tS^\perp$ contains a vector~$v$ of square~$2$, then
$\CH$ is symmetric.
\endproposition

\proof
The reflection~$t_v$
reverses the orientation of one and, hence, any maximal
positive definite subspace. On the other hand, it is obviously
an automorphism of~$\CH$, as it acts identically on~$\tS$.
\endproof

If a lattice~$N$ is unique in its genus, the existence of a vector
$v\in N$ of square~$2$ can easily be expressed in terms of
discriminant forms. Indeed, either one has
$\<v>\oplus\<v>^\perp=N$ or $\<v>\oplus\<v>^\perp\subset N$ is a
sublattice of index~$2$. In both cases, the discriminant
$\discr\<v>^\perp$ is determined by that of~$N$, and the question
reduces to the existence of a lattice~$\<v>^\perp$ within a
prescribed genus, see Theorem~\ref{N.existence}. If it does exist,
Proposition~\ref{extension.unimodular} implies that
$\<v>\oplus\<v>^\perp$ is a sublattice of a lattice isomorphic
to~$N$. Next statement is a simple special case of the above
observation.

\proposition\label{existence.-2}
Any indefinite even lattice~$N$ with
$\rank N\ge\ell(\CN)+2$ has a vector of square~$2$.
\endproposition

\proof
First, let $\rank N\ge3$. Then Theorem~\ref{N.uniqueness} implies
that $N$ is unique in its genus. On the other hand, from
Theorem~\ref{N.existence} it follows that there exists a lattice
$N'$ of signature $(\Gs_+N-1,\Gs_-N)$ whose discriminant form is
$\CN\oplus\<-\frac12>$.
(Since both signature and
Brown invariant are additive, condition~\iref{N.existence}{Br}
holds automatically.) Then the sum $\<v>\oplus N'$, $v^2=2$, is an
index~$2$ sublattice in a lattice isomorphic to~$N$.

In the exceptional case $\rank N=2$ one has $\ell(\CN)=0$, \ie,
$N$ is unimodular. Then $N\cong\bU$ and the statement is obvious.
\endproof

\subsection{Sextics with few singularities}\label{S.few.points}
Let $\mu=\rank\SS=\rank\tS-1$ be the total Milnor number of the
singularities. One has
$\mu\le19$, and the orthogonal complement~$\tS^\perp$ has
rank $21-\mu$ and signature $(2,19-\mu)$.

\theorem\label{th.l+mu<=19}
Each configuration $\tS\supset S=\SS\oplus\<h>$ satisfying the
inequality $\ell(\tCS)+\rank\SS\le19$ is realized by a unique
rigid isotopy class of plane sextics.
\endtheorem

\proof
The inequality $\ell(\tCS)+\rank\SS\le19$
implies that $\mu=\rank\SS\le18$, as otherwise $\tS^\perp$ would
be a unimodular even lattice of signature~$(2,0)$.
Thus, $\tS^\perp$ is indefinite,
$\rank\tS^\perp\ge3$, and $\rank\tS^\perp\ge\ell(\discr\tS^\perp)+2$.
Hence, Theorem~\ref{N.strong.uniqueness}
applies to~$\tS^\perp$ and $\tS$ extends to a unique abstract
homological type,
which is symmetric due to
Propositions~\ref{symmetric.-2} and~\ref{existence.-2}.
(The
existence of a lattice~$\tS^\perp$ realizing the given genus
follows from Theorem~\ref{N.existence}, as
condition~\iref{N.existence}{Br} holds automatically.)
\endproof

\corollary\label{cor.l+mu<=19}
Each configuration extending a set of singularities~$\SS$
satisfying the inequality $\ell(\discr\SS)+\rank\SS\le19$ is
realized by a unique rigid isotopy class of plane sextics.
\endcorollary

\proof
Let~$\tS$ be a configuration extending~$\SS$. Then, for each
prime~$p\ne2$, one has
$\ell_p(\tCS)\le\ell_p(\discr\SS)\le19-\mu$. For $p=2$ one has
$\ell_2(\tCS)\le\ell_2(\discr\SS)+1\le20-\mu$. However, since
$\ell_2(\CN)=\rank N\bmod2$ for each lattice~$N$, the latter
inequality still implies $\ell_2(\tCS)\le19-\mu$, and
Theorem~\ref{th.l+mu<=19} applies.
\endproof

\subsection{Classical Zariski pairs}\label{S.Zariski}
Consider a set of singularities of a classical Zariski pair of
irreducible curves, \ie, a set~$\SS$ as in~\eqref{Zariski.set}.
Let
$$
g(\SS)=10-3e-\sum_{i=1}^6 a_i\Bigl[\frac{3i}2\Bigr]-n
$$
be its virtual genus.
In~\cite{poly} it is
conjectured that,
if $g(\SS)>0$ (respectively, $g(\SS)=0$), then~$\SS$
is realized by exactly two (respectively, one) rigid isotopy
classes of irreducible sextics, one abundant and one not
(respectively, one abundant). We prove the conjecture in the case
$n=0$.

\Remark
Now, it seems clear that the nonexistence part of the
conjecture (the case $g(\SS)=0$) is wrong: a simple estimate and
Theorem~\ref{N.existence} show that \emph{most} sets of
singularities
are realized by both abundant and non-abundant curves.

The uniqueness part seems to follow more or less directly from
Theorem \ref{N.strong.uniqueness} for all abundant curves, as well
as for all curves with $e\le1$. However, as there still are quite
a number of details to be double checked (and the non-abundant
case with $e>1$ requires tedious manual calculations, \cf.
the proof of Theorem~\ref{th.Zariski}), I will consider the
general case
in a separate paper.
\endRemark

\theorem\label{th.Zariski}
Any set of singularities~$\SS$ of the form
$$
\SS=e\bE_6\oplus\bigoplus_{i=1}^6a_i\bA_{3i-1},\quad
2e+\sum ia_i=6,
$$
is realized by exactly two rigid isotopy classes of irreducible
plane sextics, one abundant and one not.
\endtheorem

\proof
We need to show that, assuming that the number of nodes $n=0$,
each of the two configurations given by
Corollary~\ref{Zariski.conf} extends to a unique homological type,
which is symmetric. One has $\mu=\rank\SS=18-(m-e)$, where $m$ is
the total number of points in~$\SS$. Since each singular point
considered
contributes exactly one to both $\ell(\discr\SS)$ and
$\ell_3(\discr\SS)$, one has $\ell(\discr\SS)+\mu=18+e$. Thus, if
$e=0$ or~$1$, the statement of the theorem follows directly from
Corollary~\ref{cor.l+mu<=19}. It
remains to consider the three sets $\SS=2\bE_6\oplus\bA_5$,
$2\bE_6\oplus2\bA_2$, or $3\bE_6$ corresponding to $e=2,3$.

If $\tS$ is an abundant configuration without points of
type~$\bA_7$ or~$\bA_{17}$, then
$\ell_3(\tS)=\ell_3(\discr\SS)-2$. If also $e\ge2$, then
still $\ell(\tCS)+\mu\le19$ and
Theorem~\ref{th.l+mu<=19} applies.

The remaining three cases, the non-abundant configurations $\tS=S$
with $e=2$ or~$3$, are considered below.

In the first two cases ($e=2$), the uniqueness of the
lattice~$\tS^\perp$ in its genus (Step~2 in
Section~\ref{S.arithmetics}) follows from
Theorem~\ref{N.uniqueness}. We will show that the homomorphism
$\OAut{\tS^\perp}\to\Aut\tCS$ is onto and that $\tS^\perp$ has a
$+$-disorienting isometry whose image in $\Aut\tCS$ belongs to
the product of $\Aut_h\tCS$ and the image of $\Aut^+\tS^\perp$
(see Steps~3 and~4 in Section~\ref{S.arithmetics}).

\subsubsection**{The case $\SS=2\bE_6\oplus\bA_5$}
One has
\botsmash{$\discr\tS^\perp\cong-\tCS\cong2\<-\frac12>\oplus3\<-\frac23>$},
so that one can take $\tS^\perp=\<-2>\oplus\<6>\oplus\bU(3)$.
Let~$a$ and~$b$ be generators of the $\<-2>$- and $\<6>$-summands,
respectively, and let $u,v$ be a standard basis for the
$\bU(3)$-summand.

Consider the $3$-primary part
$\discr\tS^\perp\otimes\Z_3\cong3\<-\frac23>$. Its automorphisms
are permutations of the three generators and
multiplication of some of them by~$(-1)$.
One can take for the generators the classes
$\beta^\pm=[b^\pm/3]$ and $\gamma=[c/3]$, where
$b^\pm=b\pm(u-2v)$ and $c=u-v$ are vectors of square~$(-6)$. Then
the reflections $\refl{b^\pm}$ and $\refl{c}$, which are well
defined elements of $\OPlus\tS^\perp$, act on $\discr\tS^\perp$ by
multiplying the corresponding generators by~$(-1)$. The reflection
$\refl{u+v}$ transposes~$\beta^+$ and~$\beta^-$, and
the reflection $\refl{b+u}$ transposes~$\beta^-$ and~$\gamma$.
Since $(u+v)^2=(b+u)^2=6$, the latter two reflections are
$+$-disorienting. All isometries mentioned act identically
on $\discr\tS^\perp\otimes\Z_2$, and together they generate the
group $\Aut(\discr\tS^\perp\otimes\Z_3)$.

The $2$-primary part is
\botsmash{$\discr\tS^\perp\otimes\Z_2\cong2\<-\frac12>$}.
Its only nontrivial automorphism is
realized by the reflection $\refl{a+b}$, which acts identically on
$\discr\tS^\perp\otimes\Z_3$.

Since the homomorphism $\OAut{\tS^\perp}\to\Aut\discr\tCS^\perp$ is
onto, one can assume that $-\discr\tCS^\perp$ and $\tCS$ are
identified so that $\beta^+$ and $\beta^-$ are generators of
the two copies of $\discr\bE_6$ in~$\tCS$. Then they can be
transposed by an admissible isometry of~$\tS$. On the other
hand,
the transposition
$\beta^+\leftrightarrow\beta^-$
is realized by a $+$-disorienting isometry
of~$\tS^\perp$. Hence, the abstract homological type is symmetric.

\subsubsection**{The case $\SS=2\bE_6\oplus2\bA_2$}
One has
\botsmash{$\discr\tS^\perp\cong-\tCS\cong\<-\frac12>\oplus4\<-\frac23>$},
so that
$\tS^\perp=\<-2>\oplus2\bU(3)$. Let~$c$ be a generator of the
$\<-2>$-summand, and let $(u_1,v_1)$ and $(u_2,v_2)$ be some
standard bases for the two $\bU(3)$-summands.

It suffices to consider the $3$-primary part
$\discr\tS^\perp\otimes\Z_3$. One can take for a basis the classes
$\alpha_i=[a_i/3]$ and $\beta^\pm=[b^\pm/3]$, where $a_i=u_i-v_i$,
$i=1,2$, and $b^\pm=(u_1-2v_1)\pm(u_2+v_2)$ are vectors of
square~$(-6)$. The reflections $\refl{a_i}$ and $\refl{b^\pm}$
multiply the corresponding generators by~$(-1)$, and modulo these
automorphisms each vector of square $(-2/3)$ in $\discr\tS^\perp$
is either one of the four generators or their sum
$\alpha_1+\alpha_2+\beta^++\beta^-=[(v_1+u_2-v_2)/3]$. Thus, each
element $\chi\in\discr\tS^\perp$ of square $(-2/3)$ can be realized
by a vector $x\in\tS^\perp$ of square~$(-6)$. The orthogonal complement
$\<x>^\perp_{\tS^\perp}$
has the genus of the lattice~$\tS^\perp$
considered in the previous case;
due to the results
obtained there,
any such element~$\chi$ can be taken to~$\alpha_1$,
and then any automorphism of the complement
$\<\alpha_1>^\perp_{\smash{\discr\tS^\perp}}$ can be realized by an
isometry of
$\<a_1>^\perp_{\tS^\perp}$.

A consideration similar to the previous case shows that
$\tS^\perp$ has a $+$-disorienting isometry that extends to an
admissible isometry of~$L$.

\subsubsection**{The case $\SS=3\bE_6$}
One has $\tCS\cong\<\frac12>\oplus3\<\frac23>$. The automorphisms
of~$\tCS$ are permutations of the three generators of order~$3$ or
multiplication of some of them by~$(-1)$. Each such automorphism
is realized by an admissible isometry of~$\tS$ (respectively,
permutation of the $\bE_6$-components and the nontrivial
symmetries of the Dynkin graphs of some of them), and it remains
to show that $\tS^\perp$ is unique in its genus and has a
$+$-disorienting isometry.
The uniqueness in the genus follows from
Theorem~\ref{N.uniqueness}; one can take
$\tS^\perp=\<6>\oplus\bU(3)$, and the
generator of the $\<6>$-summand
defines a $+$-disorienting reflection.
\endproof

\Remark
The case of abundant unnodal curves is treated in~\cite{poly}
geometrically. Since unnodal curves of the form $f_2^3+f_3^2=0$
are generic (for given
polynomials~$f_2$, $f_3$), this case reduces to the classification
of certain reducible quintics. This is done in~\cite{quintics}.
Thus, essentially new in Theorem~\ref{th.Zariski} is the case of
non-abundant curves.
\endRemark

\subsection{Maximal sextics}\label{S.examples}
The maximal value of the total Milnor number $\mu=\rank\SS$
is~$19$. When $\mu=19$, the orthogonal
complement~$\tS^\perp$ is a positive definite lattice of rank~$2$.
This case, although not covered by general
Theorems \ref{N.uniqueness} and~\ref{N.strong.uniqueness}, can
easily be handled directly, see Section~\ref{S.rank2}.
It is rather
straightforward to extend Yang's algorithm~\cite{Yang} listing the
maximal sets of singularities and corresponding configurations in
order to produce a listing of all maximal (in the sense $\mu=19$)
rigid isotopy classes. I am planning to publish the results in a
forthcoming paper.

Below, we treat manually a few examples. I have chosen the sets of
singularities that were studied in details in a recent series of
papers by Artal \all., see~\cite{Artal.first}--\cite{Artal.last},
so that the classification obtained arithmetically could be
compared with the known geometric properties of the curves. The
principal purpose of this section is to illustrate the phenomena
that take place and the number of details that should be taken
into account in an attempt to realize the algorithm
programmatically.

In the proofs below, the isomorphism classes within a given genus
for~$\tS^\perp$ can be found \via\ {\tt Maple},
using Section~\ref{S.rank2}.
Indeed, any lattice $N=\bM(a,b,c)$ with a given discriminant
form~$\CN$ must satisfy
$\mathopen|\CN\mathclose|/4\le ac\le\mathopen|\CN\mathclose|/3$,
and it remains to enumerate the triples $(a,b,c)$, calculate the
discriminant forms, and compare them to~$\CN$.
We merely indicate the result by a
sentence like `the only possibility for~$\tS^\perp$ is \dots.'

\proposition\label{D19}\label{examples.first}
The set of singularities $\SS=\bD_{19}$ extends to a unique
abstract homological type, which is symmetric.
\endproposition

\proof
One has $\CS\cong\<-\frac34>\oplus\<\frac12>$. This form has no
isotropic subgroups; hence, always $\tS=S$, and the only
possibility for~$\tS^\perp$ is $\bM(1,0,2)$.
The only nontrivial automorphism of
$\discr\tS^\perp$ is the multiplication of an order~$4$ element
by~$(-1)$; it is realized by a reflection in~$\tS^\perp$. Hence,
there is a unique abstract homological type, and it is symmetric
due to Proposition~\ref{symmetric.-2}.
\endproof

\proposition\label{A19}
The set of singularities $\SS=\bA_{19}$ admits a unique
configuration~$\tS$. It extends to two
abstract homological types, which are both
symmetric. The two lattices $\tS^\perp$ are isomorphic.
\endproposition

\Remark
Sextics with this set of singularities were studied in Artal
\all.~\cite{Artal.Trends}, where the two rigid isotopy classes
were discovered. The only difference between the two
homological types is the anti-isometry identifying the
discriminant groups of~$\tS$ and~$\tS^\perp$. (One also observes a
similar phenomenon in Propositions~\ref{2A9+A1}
and~\ref{E6+A7+A3+A2+A1} below, where the curves are reducible.)
Thus, it \emph{looks} like the two pairs in question are obtained
by gluing diffeomorphic pieces \via\ different diffeomorphisms of
their boundaries, \cf. Remark~\ref{rem.M-V}. At present, I do
not know
whether this claim is true, as of course the global
Torelli theorem only applies to whole $K3$-surfaces.

Up to projective equivalence, each rigid isotopy class consists of
a single curve~$C_i$, $i=1,2$, and the two curves are indeed very
similar to each other. For example, disregarding the hyperplane section
class~$h$ in the calculation above, one can easily  show that the
two covering $K3$-surfaces~$\BX_i$ are deformation equivalent. The
fundamental groups $\pi_1(\Cp2\sminus C_i)$ were calculated
in~\cite{Artal.Trends}, and it was shown that they are isomorphic
to each other. Whether the two complements $\Cp2\sminus C_i$
themselves are diffeo-/homeomorphic still remains an open question.
\endRemark

\proof
One has
$\CS\cong\<-\frac{19}{20}>\oplus\<\frac12>\cong
 \<\frac45>\oplus\<\frac14>\oplus\<\frac12>$,
the first two
summands being generated by~$4\alpha$ and~$5\alpha$, where
$\alpha$ is a canonical generator of the group
$\discr\bA_{19}$. Since $\CS$ has no isotropic subgroups, one has
$\tS=S$, and the only possibility for the orthogonal
complement~$\tS^\perp$ is $\bM(1,0,10)$. The automorphism group
$\Aut\CS\cong\Z/2\Z\times\Z/2\Z$
consists of the automorphisms
$$
(4\alpha,5\alpha)\mapsto(\epsilon_1\cdot4\alpha,\epsilon_2\cdot5\alpha),\quad
 \epsilon_1,\epsilon_2=\pm1,
$$
whereas the images of both $\OAut\SS$ and $\OAut{\tS^\perp}$ are
generated by $-\id$, corresponding to
$\epsilon_1=\epsilon_2=-1$. Hence, $\tS$ extends to two distinct
abstract homological types. They are both symmetric due to
Proposition~\ref{symmetric.-2}.
\endproof

\proposition\label{A18+A1}
The set of singularities $\SS=\bA_{18}\oplus\bA_1$ admits a unique
configuration~$\tS$. It extends to two
abstract homological types, which differ by the lattice
$\tS^\perp$. One of the homological types is symmetric,
the other one is not,
so that there are three rigid
isotopy classes of sextics with this set of singularities.
\endproposition

\Remark
This set of singularities was first studied in
Artal \all.~\cite{Artal.Trends}. 
The most remarkable fact is the existence of two rigid isotopy
classes that
differ by the orientation of their homological types.
This example may be a first candidate for a pair of
sextic curves~$C_1$, $C_2$ with homeomorphic but not diffeomorphic
pairs $(\Cp2,C_i)$,
see Remark~\ref{rem.diffeo}.
According to the description of the orthogonal
groups given in Section~\ref{S.rank2}, this situation should be rather
typical for the maximal Milnor number: any abstract
homological type with $\tS^\perp\cong\bM(a,b,c)$, $0<b<a<c$, would
be asymmetric.

All three curves are given by Galois conjugate equations defined
over $\Q(\beta)$, where $\beta$ is a root of $19s^3+50s^2+36s+8$.
In~\cite{Artal.Trends}, there are more examples of curves given by
complex conjugate equations; the corresponding sets of
singularities are
$\bA_{16}\oplus\bA_2\oplus\bA_1$ and $\bA_{15}\oplus\bA_4$.
\endRemark

\proof
One has
$\CS\cong\<-\frac{18}{19}>\oplus\<-\frac12>\oplus\<\frac12>$. The
only imprimitive extension of~$S$ would contradict
Proposition~\ref{configuration}; hence, there is a unique
configuration \topsmash{$\tS=S$}. The genus of~$\tS^\perp$ contains two
isomorphism classes: $\bM(1,0,19)$ and $\bM(4,2,5)$. The only
nontrivial automorphism of~$\tCS$ is realized by the isometry
$-\id\in\OAut{\tS^\perp}$. Hence, each of the isomorphism classes
gives rise to a unique abstract homological type. The abstract
homological type
with $\tS^\perp=\bM(1,0,19)$ is symmetric due to
Proposition~\ref{symmetric.-2}; the one
with $\tS^\perp=\bM(4,2,5)$
is not symmetric as $\tS^\perp$ has no
$+$-disorienting isometries.
\endproof

\proposition\label{2A9+A1}
The set of singularities $\SS=2\bA_9\oplus\bA_1$ admits two
distinct configuration~$\tS$, with $[\tS\mathbin:S]=2$ or~$10$.
The former configuration extends to two
abstract homological types \rom(with isomorphic
lattices~$\tS^\perp$\rom), the latter extends to one. All
extensions are symmetric, so that altogether there are three rigid
isotopy classes of sextics with this set of singularities.
\endproposition

\Remark
Although the configurations differ, all three curves have the same
combinatorial data. The case $[\tS\mathbin:S]=10$ can be told
apart by the existence of an extra line in a special position with
respect to the curve, see Artal \all.~\cite{Artal.KT}.
\endRemark

\proof
One has
$\CS\cong2\<-\frac9{10}>\oplus\<-\frac12>\oplus\<\frac12>\cong
 2\<\frac25>\oplus3\<-\frac12>\oplus\<\frac12>$.
Let~$\alpha_{1,2}$, $\beta$, and~$\gamma$ be generators of the
summands $\discr\bA_9$, $\discr\bA_1$, and $\discr\<2>$,
respectively.

Since $\ell_2(\CS)=4$, the kernel
\topsmash{$\CK=\tS/S$} must contain elements of order~$2$. The
isotropic elements of order~$2$ are $\beta+\gamma$ and
$5\alpha_{1,2}+\gamma$. The former cannot belong to~$\CK$ due to
Proposition~\ref{configuration};
the two latter are interchangeable by an admissible isometry.
Hence, one can assume that $\CK\otimes\Z_2\cong\Z/2\Z$ is
generated by $5\alpha_1+\gamma$.

The $5$-primary part $\CK\otimes\Z_5$ may be either trivial or one
of the two order~$5$ subgroups generated by
$2\alpha_1\pm4\alpha_2$. In the latter case the two subgroups are
conjugate by an admissible isometry preserving
$5\alpha_1+\gamma$. Hence, up to admissible isometry, there
are two configurations, which differ by the index
$[\tS\mathbin:S]$.

If $[\tS\mathbin:S]=10$, then $\tCS\cong2\<-\frac12>$ and
$\tS^\perp\cong\bM(1,0,1)$. The homomorphism
\topsmash{$\OAut{\tS^\perp}\to\Aut\discr\tS^\perp$} is onto, and the
only resulting abstract homological type is symmetric due to
Proposition~\ref{symmetric.-2}.

If $[\tS\mathbin:S]=10$, then
$\tCS\cong2\<\frac25>\oplus2\<-\frac12>$ and
$\tS^\perp\cong\bM(5,0,5)$. The group
\topsmash{$\Aut\tCS$} is generated by the multiplications
$(-1)_i\:\alpha_i\mapsto-\alpha_i$, $i=1,2$, the
transposition~$\tr_5$ of~$\alpha_1$ and~$\alpha_2$, and the
transposition~$\tr_2$ of the generators of the two summands of
order~$2$. The subgroup $\Aut_h\tCS$ is generated by~$(-1)_1$
and~$(-1)_2$, and the image of~$\OAut{\tS^\perp}$ in $\Aut\tCS$ is
generated by~$(-1)_1$, $(-1)_2$, and
the composition $\tr_5\circ\tr_2$. Hence,
$\tS$ extends to two distinct abstract homological types, and they
are both symmetric (as $\tS^\perp$ has a $+$-disorienting
isometry that extends to~$L$ by an admissible isometry
of~$\tS$).
\endproof

\proposition\label{E6+A7+A3+A2+A1}\label{examples.last}
The set of singularities
$\SS=\bE_6\oplus\bA_7\oplus\bA_3\oplus\bA_2\oplus\bA_1$ admits a
unique configuration~$\tS$. It extends to two abstract homological
types, which are both symmetric. The two lattices $\tS^\perp$ are
isomorphic.
\endproposition

\Remark
Sextics with this set of singularities are all reducible,
splitting to a singular quintic and a line. They
were studied in Artal \all.~\cite{Artal.braids}.
\endRemark

\proof
One has
$\CS\cong\<\frac23>\oplus\<-\frac78>\oplus\<-\frac34>
 \oplus\<-\frac23>\oplus\<-\frac12>\oplus\<\frac12>$.
Since
$\ell_2(\CS)=4$, the kernel~$\CK$ must contain elements of
order~$2$. The only isotropic element of order~$2$ not
contradicting
Proposition~\ref{configuration}
is $4\beta_7+\beta_1+\gamma$, where
$\beta_i$ is a generator of $\discr\bA_i$, $i=1,2,3,7$, and
$\gamma$ is the generator of~$\discr\<2>$. Thus, there is a unique
configuration~\topsmash{$\tS$}, and the group
$\tCS\cong\<\frac23>\oplus\<-\frac38>\oplus\<-\frac34>\oplus\<-\frac23>$
is generated by~$\alpha$ (a generator of $\discr\bE_6$),
$\beta_7'=\beta_7+\gamma$, $\beta_3$, and~$\beta_2$. The group
\topsmash{$\Aut\tCS$} is generated by the multiplications of generators
by~$(-1)$, which lift to admissible isometries of~$S$, and the
involution
$\Gf\:(\beta_7',\beta_3)\mapsto(3\beta_7'+2\beta_3,\beta_3+4\beta_7')$,
which is not in the image of~$\Oh S$.

The only possibility for~$\tS^\perp$ is
$\bM(6,0,12)$. Since the
involution~$\Gf$ above does not lift to ~$\tS^\perp$,
there are two abstract homological types extending~$\tS$. Each of
them is symmetric, as the $+$-disorienting reflection defined by
the vector of square~$24$ extends to~$L$ by an admissible
isometry of~$\tS$.
\endproof

\subsection{Concluding remarks}\label{S.remarks}

\subsubsection{Examples with $\mu<19$}\label{Q.mu<19}
One of the by-products of the calculation in
Section~\ref{S.examples} is the fact that each of the four steps
outlined in Section~\ref{S.arithmetics} does matter, in the sense
that there are pairs of sextics that diverge at that particular
step. However, in all these examples one has $\mu=19$, \ie, the
moduli space is discrete, and I do not know a single
examples of a pair of not rigidly isotopic sextics with $\mu<19$
that share the same configuration. (A number of families with
$\mu=18$ is considered in~\cite{Artal.Trends}, where it is
proved that the curves are determined by their combinatorial data.
The members of classical Zariski pairs considered in
Section~\ref{S.Zariski} differ by their configurations.)

At present, I do not know how general this phenomenon is and how
it could be proved/disproved without essentially enumerating all
rigid isotopy classes.

\subsubsection{Asymmetric homological types}\label{Q.asymmetric}
The existence of two opposite orientations of periods of marked
$K3$-surfaces is a well known fact. The two orientations are
interchanged by the canonical real structure on the moduli space,
sending a $K3$-surface to its conjugate. Thus, asymmetric abstract
homological types give rise to moduli spaces without real points.
It would be interesting to find similar examples with $\mu<19$,
when the modular space has positive dimension, or to prove that such
examples do not exist, \cf. Section \ref{Q.mu<19}. (Note that,
\emph{a priori}, the moduli space may have no real points even if
the homological type is symmetric.)

For plane curves, the canonical real structure on the moduli space is
induced by the standard complex conjugation on~$\Cp2$. In
particular, the pairs $(\Cp2,C_i)$ corresponding to two conjugate
curves \emph{are} diffeomorphic, but the diffeomorphism is not
regular and, most importantly, it reverses the orientation of
complex curves, inducing~$(-1)$ in $H_2(\Cp2)$, \cf.
Section~\ref{Q.regularity}.

Apparently, the existence of asymmetric homological types is due
to the fact that we consider $K3$-surfaces with a fixed
polarization (which prohibits certain obvious changes of the
marking), whereas no involution that would interchange the
components is assumed \emph{a priori} (as in the case of real
curves and surfaces).
Probably, one can
anticipate a similar phenomenon in the case of quartic surfaces
in~$\Cp3$, see Section~\ref{Q.quartics}.

\subsubsection{Conjugacy over number fields}\label{Q.conjugacy}
Comparing the results obtained in Section~\ref{S.examples} with
the geometric properties of the curves discovered
in~\cite{Artal.first}--\cite{Artal.last}, one can observe that, in
the case of maximal total Milnor number $\mu=19$, {\proclaimfont
all sextics with a given configuration~$\tS$ are given by
equations Galois conjugate over a certain finite extension
of~$\Q$.} (Remind that curves with $\mu=19$ are rigid, \ie, their
moduli spaces are discrete.) At present, I do not know how general
this statement is and how it can be obtained arithmetically.
Another piece of substantiating evidence is given by the material
of~\S\ref{SS.sextics}, where various geometric properties of curves that
should remain invariant under all, not necessarily continuous,
Galois transformations are expressed in terms of the configuration
only.

\subsubsection{The regularity condition}\label{Q.regularity}
The regularity condition in Theorem \ref{th.main} seems to be a
purely technical assumption; it is needed to assure that the
diffeomorphism~$f$ lifts to the minimal resolutions of
singularities and then, further, to the nonsingular double
coverings. I believe that any diffeomorphism preserving the
complex orientations gives rise to a regular one. To
prove/disprove this statement one would need to know the
diffeotopy classification of auto-diffeomorphisms of the link
(\ie, boundary of a regular neighborhood) of each simple
singularity of surfaces. Unfortunately, I do not know any results
in this direction.

Certainly, if it turns out that the regularity assumption can be
dropped, one would still have to require that $f$ preserve the
complex orientation of both~$\Cp2$ and the curves, \ie, that the
induced homomorphism $f_*\:H_*(\Cp2)\to H_*(\Cp2)$ is the
identity.

\subsubsection{Other equivalence relations}\label{Q.other.relations}
In the definition of admissible iso\-metry, it is required that the
distinguished basis~$\Gs$ of~$\SS$ should only be fixed as a
set. Geometrically, this means that neither the order of the
singular points nor their `orientation' are assumed fixed. If the
order matters, one should fix the basis and consider isomorphism
classes of lattice polarized $K3$-surfaces in the sense of
Nikulin~\cite{NikK3}; they are classified by oriented abstract
homological types up to isometries identical on $\SS\oplus\<h>$.

An interesting example of multiple equivalence classes is
described in Remark \ref{rem.Zariski.orientation}. A more
straightforward example is given by Proposition~\ref{2A9+A1},
where the two $\bA_9$ points cannot be transposed by a rigid
isotopy: they are distinguished by the combinatorial data of the
curves.

Alternatively, one can compare curves using various relaxed
equivalence relations: homeomorphism of the pairs,
diffeo-/home\-o\-mor\-phism of the complement spaces, \etc. For
most curves, this problem still remains open.

\subsubsection{Quartic surfaces}\label{Q.quartics}
The classification of singular quartic surfaces in~$\Cp3$ should
be very similar to the classification of plane sextics. The proof
of the corresponding counterpart of Theorem~\ref{th.main} would
repeat literally the contents of \S\ref{SS.moduli}, with $h^2=2$
replaced with $h^2=4$. It is worth mentioning that, as in the case
of plane sextics, the rigid isotopy class of a quartic surface
with at least one non-simple singular point is determined by its
combinatorial data, see~\cite{quartics2}.

\refstyle{C}

\enddocument